\documentclass[12pt,epsf,bezier]{article}
\renewcommand{\baselinestretch}{1.7}
\usepackage{latexsym}
\usepackage{epsfig}

\begin{document}
\newtheorem{pro}{Proposition}
\newtheorem{theo}{Theorem}
\newtheorem{coro}{Corollary}
\newtheorem{lem}{Lemma}

\newcommand{\bboldy}[1]{\mbox{\boldmath${\#1}$}}

\def\E{\mbox{\rm E}}
\def\Var{\mbox{\rm Var}}
\def\Cov{\mbox{\rm Cov}}
\def\Corr{\mbox{\rm Corr}}
\def\Pr{\mbox{\rm Pr}}
\def\I{\mbox{\rm I}}
\def\bbeta{\mbox{\boldmath${\beta}$}}
\def\bgamma{\mbox{\boldmath${\gamma}$}}
\def\sbeta{\mbox{\scriptsize \boldmath${\beta}$}}
\def\eeps{\mbox{\boldmath${\epsilon}$}}
\def\bmu{\mbox{\boldmath${\mu}$}}
\def\bxi{\mbox{\boldmath${\xi}$}}
\def\real{\hbox{\rm\setbox1=\hbox{I}\copy1\kern-.45\wd1 R}}
\newcommand{\bgams}{\mbox{\boldmath{\scriptsize $\gamma$}}}
\newcommand{\bLam}{\mbox{\boldmath{$\Lambda$}}}
\newcommand{\Ell}{\mbox{$\cal{L}$}}
\newcommand{\cG}{\mathcal{G}}
\newcommand{\cL}{\mathcal{L}}
\newcommand{\half}{\frac{1}{2}}
\newcommand{\bsh}{\parindent 0em}
\newcommand{\esh}{\parindent 2.0em}
\newcommand{\eps}{\epsilon}
\newcommand{\peps}{{(\eps)}}
\def\bU{{\bf U}}
\def\bV{{\bf V}}
\def\bG{{\bf G}}
\def\bC{{\bf C}}
\def\bD{{\bf D}}
\def\bu{{\bf u}}

\begin{center}
\bf
\large
Prospective survival analysis with a general
semiparametric shared frailty model - a pseudo full likelihood
approach
\end{center}

\normalsize

\vspace*{1cm}

\begin{center}

\bf

Malka Gorfine\footnote{ To whom correspondence should be
addressed} \\ {\it Department of Mathematics, Bar-Ilan University,
Ramat-Gan 52900, Israel} \\ {gorfinm@macs.biu.ac.il}
\\

\vspace*{1cm}
David M. Zucker\\
{\it Department of Statistics, Hebrew University,
Mt. Scopus, Jerusalem 91905, Israel}\\
{mszucker@mscc.huji.ac.il} \\

\vspace*{1cm}
Li Hsu \\
{\it Division of Public Health Sciences, Fred Hutchinson
Cancer Research Center, Seattle, WA 98109-1024 USA}\\
{lih@fhcrc.org}
\end{center}

\vspace*{1in}

\begin{center}
\today
\end{center}

\newpage
\section*{Summary}
In this work we provide a simple estimation procedure for a
general frailty model for analysis of prospective correlated
failure times. Rigorous large-sample theory for the proposed
estimators of both the regression coefficient vector and the
dependence parameter is given, including consistent variance
estimators. In a simulation study under the widely used gamma
frailty model, our proposed approach was found to have essentially
the same efficiency as the EM-based estimator considered by other
authors, with negligible difference between the standard errors of
the two estimators. The proposed approach, however, provides a
framework capable of handling general frailty distributions with
finite moments and yields an explicit consistent variance
estimator.\\

\noindent {\it Key words:}  Correlated failure times; EM
algorithm; Frailty model; Prospective family study; Survival
analysis.

\newpage
\section{Introduction}
Many epidemiological studies involves failure times that are
clustered into groups, such as families or schools, where some
unobserved characteristics shared by the members of the same
cluster (e.g. genetic information or unmeasured shared
environmental exposures) could influence time to the studied
event. In frailty models within cluster dependence is represented
through a shared unobservable variable as a random effect.
Estimation in the frailty model has received much attention under
various frailty distributions, including gamma (Gill, 1985, 1989;
Nielsen et al., 1992; Klein 1992, among others), positive stable
(Hougaard, 1986; Fine et al., 2003), inverse Gaussian, compound
Poisson (Henderson and Oman, 1999) and log-normal (McGilchrist,
1993; Ripatti and Palmgren, 2000; Vaida and Xu, 2000, among
others). Hougaard (2000) provides a comprehensive review of the
properties of the various frailty distributions. In a frailty
model, the parameters of interest typically are the regression
coefficients, the cumulative baseline hazard function, and the
dependence parameters in the random effect distribution.

Since the frailties are latent covariates, the
Expectation-Maximization (EM) algorithm is a natural estimation
tool, with the latent covariates estimated in the E-step and the
likelihood maximized in the M-step by substituting the estimated
latent quantities. Gill (1985), Nielsen et al. (1992) and Klein
(1992) discussed EM-based maximum likelihood estimation for the
semiparametric gamma frailty model. One problem with the EM
algorithm is that variance estimates of the estimated parameters
are not readily available (Louis, 1982; Gill, 1989; Nielsen et
al., 1992; Andersen et al., 1997). It was suggested (Gill, 1989;
Nielsen et al, 1992) that a nonparametric information calculation
could yield consistent variance estimators. Parner (1998),
building on Murphy (1994, 1995), proved the consistency and
asymptotic normality of the maximum likelihood estimator in the
gamma frailty model. Parner also presented a consistent estimator
of the limiting covariance matrix of the estimator based on
inverting a discrete observed information matrix. He noted that
since the dimension of the observed information matrix is the
dimension of the regression coefficient vector plus the number of
observed survival times, inverting the matrix is practically
infeasible for a large data set with many distinct failure times.
Thus, he proposed another covariance estimator based on solving a
discrete version of a second order Sturm-Liouville equation. This
covariance estimator requires substantially less computational
effort, but still is not so simple to implement.

The purpose of our work here is to develop a new inference
technique that can handle any parametric frailty distribution with
finite moments. Our new method possesses a number of desirable
properties: a non-iterative procedure for estimating the
cumulative hazard function; consistency and asymptotic normality
of parameter estimates; a direct consistent covariance estimator;
and easy computation and implementation. The rest of the paper is
organized as follows. In Section 2 we present the estimation
procedure. Consistency and asymptotic results for the estimators
are given in Section 3. As the frailty model is often applied
using a gamma frailty distribution, Section 4 compares the finite
sample performance of our approach and the EM-based approach under
the gamma distribution. Section 5 provides an example using a
diabetic retinopathy data set. Section 6 presents concluding
remarks.

\section{The Proposed Approach}
Consider $n$ families, with family $i$ containing $m_i$ members,
$i=1,\ldots,n$. Let $\delta_{ij}=I(T^0_{ij} \leq C_{ij})$ be a
failure indicator where $T^0_{ij}$ and $C_{ij}$ are the failure
and censoring times, respectively, for individual $ij$. Also let
$T_{ij}=\min(T^0_{ij},C_{ij})$ be the observed follow-up time and
${\bf Z}_{ij}$ be a $p \times 1$ vector of covariates. In
addition, we associate with family $i$ an unobservable
family-level covariate $W_i$, the ``frailty", which induces
dependence among family members. The conditional hazard function
for individual $ij$ conditional on the family frailty $W_i$, is
assumed to take the form
 $$
\lambda_{ij}(t)=W_i \lambda_0(t) \exp(\bbeta^T {\bf Z}_{ij})
\;\;\;\;\;\; i=1,\ldots,n \;\;\; j=1,\ldots,m_i
 $$
where $\lambda_0$ is an unspecified conditional baseline hazard
and $\bbeta$ is a $p \times 1$ vector of unknown regression
coefficients. This is an extension of the Cox (1972) proportional
hazards model, with the hazard function for an individual in
family $i$ multiplied by $W_i$. We assume that, given ${\bf
Z}_{ij}$ and $W_i$, the censoring is independent and
noninformative for $W_i$ and $(\bbeta,\Lambda_0)$ (in the sense of
Andersen et al., 1993, Sec. III.2.3). We assume further that the
frailty $W_i$ is independent of ${\bf Z}_{ij}$ and has a density
$f(w;\theta)$, where $\theta$ is an unknown parameter. For
simplicity we assume that $\theta$ is a scalar, but the
development extends readily to the case where $\theta$ is a
vector. Let $\tau$ be the end of the observation period. The full
likelihood of the data then can be written as
\begin{eqnarray}\label{eq:like}
\lefteqn{L}
 &=&\Pi_{i=1}^n \int \Pi_{j=1}^{m_i} \{ \lambda_{ij}(T_{ij}) \}^{\delta_{ij}}
 S_{ij}(T_{ij}) f(w) dw \nonumber \\
 &=& \Pi_{i=1}^n \Pi_{j=1}^{m_i} \{\lambda_0(T_{ij})\exp(\bbeta^T {\bf Z}_{ij})
 \}^{\delta_{ij}} \Pi_{i=1}^n
 \int w^{N_{i.}(\tau)} \exp \{-w H_{i.}(\tau) \}f(w)dw,
\end{eqnarray}
where $N_{ij}(t)=\delta_{ij}I(T_{ij}\leq t)$,
$N_{i.}(t)=\sum_{j=1}^{m_i}N_{ij}(t)$,
$H_{ij}(t)=\Lambda_0(T_{ij}\wedge t)\exp(\bbeta^T {\bf Z}_{ij})$,
$a \wedge b = \min\{a,b\}$, $\Lambda_0(\cdot)$ is the baseline
cumulative hazard function, $S_{ij}(\cdot)$ is the conditional
survival function of subject $ij$, and
$H_{i.}(t)=\sum_{j=1}^{m_i}H_{ij}(t)$. The log-likelihood is given
by
 $$
l=\sum_{i=1}^{n} \sum_{j=1}^{m_i} \delta_{ij} \log\{
\lambda_0(T_{ij})\exp(\bbeta^T {\bf Z}_{ij}) \} + \sum_{i=1}^{n}
\log \left\{ \int w^{N_{i.}(\tau)} \exp \{-w H_{i.}(\tau)\}f(w)dw
\right\}.
 $$
The normalized scores (log-likelihood derivatives) for
$(\beta_1,\ldots,\beta_p)$ are given by
\begin{equation}
 U_r=\frac{1}{n} \sum_{i=1}^{n} \sum_{j=1}^{m_i} \delta_{ij}
 Z_{ijr} - \frac{1}{n} \sum_{i=1}^n
 \frac{\left[\sum_{j=1}^{m_i} H_{ij}(T_{ij}) Z_{ijr} \right] \int w^{N_{i.}(\tau)+1}
 \exp \{-w H_{i.}(\tau)\}f(w)dw}
 {\int w^{N_{i.}(\tau)} \exp \{ -w H_{i.}(\tau)\}f(w)dw}
\label{score}
\end{equation}
for $r=1,\ldots,p$. The normalized score for $\theta$ is
 $$
 U_{p+1}=\frac{1}{n}\sum_{i=1}^{n}
 \frac{\int w^{N_{i.}(\tau)} \exp \{ -w H_{i.}(\tau)\}f'(w)dw}
 {\int w^{N_{i.}(\tau)} \exp \{ -w H_{i.}(\tau)\}f(w)dw}
 $$
where $f'(w)=\frac{d}{d\theta}f(w)$. Let
$\bgamma=(\bbeta^T,\theta)$ and ${\bf U}
(\bgamma,\Lambda_0)=(U_1,\ldots,U_p,U_{p+1})^T$. To obtain
estimators $\hat{\bbeta}$ and $\hat{\theta}$, we propose to
substitute an estimator of $\Lambda_0$, denoted by
$\hat{\Lambda}_0$, into the equations ${\bf U}
(\bgamma,\Lambda_0)=0$.

Let $Y_{ij}(t)=I(T_{ij} \geq t)$ and let ${\mathcal F}_t$ denote
the entire observed history up to time $t$, that is
 $$
{\mathcal F}_{t} = \sigma\{N_{ij}(u), Y_{ij}(u), {\bf Z}_{ij},
i=1,\ldots,n; j=1,\ldots,m_i;0 \leq u \leq t \}.
 $$
Then, as discussed by Gill (1992) and Parner (1998), the
stochastic intensity process for $N_{ij}(t)$ with respect to
${\mathcal F}_{t}$ is given by
\begin{equation}\label{eq:inten}
\lambda_0(t)\exp(\bbeta^T {\bf Z}_{ij}) Y_{ij}(t)
\psi_{i}(\bgamma,\Lambda_0,t-),
\end{equation}
where
 $$
\psi_{i}(\bgamma,\Lambda_0,t) = \E(W_i|{\mathcal F}_{t}).
 $$
Using a Bayes theorem argument and the joint density
(\ref{eq:like}) with observation time restricted to $[0,t)$, we
obtain
 $$
 \psi_{i}(\bgamma,\Lambda,t)=\phi_{2i}(\bgamma,\Lambda,t)/
 \phi_{1i}(\bgamma,\Lambda,t),
 $$
where
 $$
 \phi_{ki}(\bgamma,\Lambda_0,t)=\int w^{N_{i.}(t)+(k-1)} \exp\{-w
 H_{i.}(t)\}f(w)dw,  \;\;\; k=1,\ldots,4.
 $$
Given the intensity model (\ref{eq:inten}), in which
$\exp(\bbeta^T {\bf Z})\psi_{i}(\bgamma,\Lambda_0,t-)$ may be
regarded as a time dependent covariate effect, a natural estimator
of $\Lambda_0$ is
 a Breslow (1974) type estimator along the lines of
 Zucker (2005). For given values of $\bbeta$ and $\theta$ we
estimate $\Lambda_0$ as a step function with jumps at the observed
failure times $\tau_k$, $k=1,\ldots,K$, with
 \begin{equation}\label{eq:lambda}
\Delta \hat{\Lambda}_0(\tau_k) =\frac
 {d_k }
 {\sum_{i=1}^{n}
 \psi_{i}(\bgamma,\hat{\Lambda}_0,\tau_{k-1})
 \sum_{j=1}^{m_i} Y_{ij}(\tau_k) \exp(\bbeta^T {\bf Z}_{ij})
 }
 \end{equation}
where $d_k$ is the number of failures at time $\tau_k$. Note that
given the intensity model (\ref{eq:inten}), the estimator of the
$k$th jump depends on $\hat{\Lambda}_0$ up to and including time
$\tau_{k-1}$. By this approach, we avoid complicating the
iterative optimization process with a further iterative scheme,
like that of Shih and Chatterjee (2002), for estimating the
cumulative hazard.

\section{Large-Sample Study}
Let $\bgamma^\circ=({\bbeta^\circ}^T,\theta^\circ)^T$ with $\bbeta^\circ$,
$\theta^\circ$ and $\Lambda_0^\circ(t)$ denoting the respective true
values of $\bbeta$, $\theta$ and $\Lambda_0(t)$, and let
$\hat{\bgamma}=({\hat{\bbeta}}^T,\hat{\theta})^T$. In Appendix A,
the conditions assumed in establishing the asymptotic properties
of $\hat{\bgamma}$ are listed and discussed.

Using arguments similar to those of Zucker (2005, Appendix A.3),
the following can be shown (see, Appendix A):
\begin{description}
\item[A.]
$\hat{\Lambda}_0(t,\bgamma)$ converges almost surely to
$\Lambda_0(t,\bgamma)$ uniformly in $t$ and $\bgamma$.
\item[B.]
${\bf U}(\bgamma,\hat{\Lambda}_0(\cdot,\bgamma))$ converges almost
surely uniformly in $t$ and $\bgamma$ to a limit ${\bf u}
(\bgamma,\Lambda_0(\cdot,\bgamma))$.
\item[C.]
There exists a unique consistent root to ${\bf
U}(\hat{\bgamma},\hat{\Lambda}_0(\cdot,\hat{\bgamma}))=0$.
\end{description}

To show that $\hat{\bgamma}$ is asymptotically normally
distributed, we write
\begin{eqnarray}
\lefteqn{0}&=&{\bf U}
(\hat{\bgamma},\hat{\Lambda}_0(\cdot,\hat{\bgamma}))\nonumber
\\ &=& {\bf U}(\bgamma^\circ,\Lambda_0^\circ) +
[{\bf U} (\bgamma^\circ,\hat{\Lambda}_0(\cdot,\bgamma^\circ))-{\bf
U} (\bgamma^\circ,\Lambda_0^\circ)] \nonumber\\ & & + [{\bf
U}(\hat{\bgamma},\hat{\Lambda}_0(\cdot,\hat{\bgamma}))- {\bf
U}(\bgamma^\circ,\hat{\Lambda}_0(\cdot,\bgamma^\circ))]. \nonumber
\end{eqnarray}
In Appendix B we analyze each of the above three terms and prove
that $n^{1/2}(\hat{\bgamma}-\bgamma^\circ)$ is asymptotically
mean-zero normally distributed, with a covariance matrix that can
be consistently estimated by the sandwich estimator
 \begin{equation}\label{eq:se}
  {\bf D}^{-1}(\hat{\bgamma})\{\hat{\bf V}(\hat{\bgamma})+\hat{\bf G}(\hat{\bgamma})+\hat{\bf C}(\hat{\bgamma})\}
  {\bf D}^{{-1}}(\hat{\bgamma})^T.
 \end{equation}
The matrix ${\bf D}$ consists of the derivatives of the $U_r$'s
with respect to the parameters $\bgamma$. ${\bf V}$ is the
asymptotic covariance matrix of ${\bf
U}(\bgamma^\circ,\Lambda_0^\circ)$, ${\bf G}$ is the asymptotic
covariance matrix of $[{\bf
U}(\bgamma^\circ,\hat{\Lambda}_0(\cdot,\bgamma^\circ))-{\bf
U}(\bgamma^\circ,\Lambda_0^\circ)]$,
 and ${\bf C}$ is the asymptotic covariance matrix between ${\bf U}(\bgamma^\circ,\Lambda_0^\circ)$
 and
$
[{\bf U}(\bgamma^\circ,\hat{\Lambda}_0(\cdot,\bgamma^\circ))-{\bf
U}(\bgamma^\circ,\Lambda_0^\circ)]$. The term ${\bf G}+{\bf C}$
reflects the added variance resulting from the need to estimate
the cumulative hazard function. All the above matrices are defined
explicitly in the Appendix.

\section{Simulation Study for the Gamma Frailty Case}
Gill (1985), Klein et al. (1992), and Nielsen et al. (1992) dealt
with the gamma frailty model by applying the EM algorithm to the
Cox partial likelihood. This methods may be interpreted as a
semi-parametric full maximum likelihood method. Murphy (1994,
1995) showed consistency and asymptotic normality for the model
without covariates, where the unknown parameters are the
integrated hazard function and the gamma frailty parameters.
Parner (1998) extended the consistency and asymptotic normality
results to the correlated gamma frailty model with covariates. In
what follows we compare our proposed method to the EM method under
the gamma frailty distribution with expectation 1 and variance
$\theta$.

The following is the EM-based estimation algorithm as given in
Nielsen et al. (1992).
\begin{description}
\item[Step I:]
Using standard Cox regression software, obtain initial estimates
of $\bbeta$ and $\Lambda_0$, taking $W_{i}=1$, $i=1,\ldots,n$
(i.e. $\theta=0$).
\item[Step II (E step):]
Using the current values of $\bbeta$, $\Lambda_0$ and $\theta$,
estimate the frailty value $W_i$ by
\begin{equation}\label{eq:em}
\hat{W}_i=\frac{N_{i.}(\tau)+\theta^{-1}}{H_{i.}(\tau)+\theta^{-1}}.
\end{equation}
\item[Step III (M step):]
Update the estimate of $\bbeta$ by fitting a Cox proportional
hazard model with covariates ${\bf Z}$ and offset term
$\log(\hat{W})$. Update the estimate of $\Lambda_0$ by the
traditional Breslow type estimator associated with the Cox model.
Update the estimate of $\theta$ by the maximum likelihood
estimator based on (\ref{eq:like}).
\item[Step IV:]
Iterate between Steps II and III until convergence.
\end{description}

Our estimation technique can be summarized by the following
algorithm.
\begin{description}
\item[Step I:]
Using standard Cox regression software, obtain initial estimates
of $\bbeta$ and as initial value for $\hat{\theta}$, let
$\hat{\theta}=0$.
\item[Step II:]
Using the current values of $\bbeta$ and $\theta$, estimate
$\Lambda_0$ using the non-iterative estimate presented by Equation
(\ref{eq:lambda}).
\item[Step III:]
Using the current estimate $\hat{\Lambda}_0$, estimate $\bbeta$
and $\theta$ by solving ${\bf U}(\bgamma,\hat{\Lambda})=0$.
\item[Step IV:]
Iterate between Steps II and III until convergence.
\end{description}

It is easy to see that under the gamma distribution for $W_i$,
\begin{equation}\label{eq:ours}
\psi_i(\bgamma,\Lambda_0,t-) =
\E(W_i|F_{t-})=\frac{N_{i.}(t-)+\theta^{-1}}{H_{i.}(t-)+\theta^{-1}}.
\end{equation}
Murphy (1994) showed that for the model without covariates, an
estimator of the cumulative hazard function based on the EM
algorithm with (\ref{eq:ours}) instead of (\ref{eq:em}) converges
to the true value of the cumulative hazard function. This result
can be extended to the case where covariates are included in the
model.

Note that in Murphy (1994), the cumulative hazard function at
$\tau_k$ includes the cumulative information up through time
$\tau_k$, whereas in the EM algorithm the accumulated information
is up through time $\tau$, the entire study period. In contrast,
in our approach the cumulative hazard function at $\tau_k$ only
includes the information up through the previous failure time
point $\tau_{k-1}$. Hence, one might suspect our estimators are
somewhat less efficient than the EM-based estimators. Part of the
goal of our simulation study was to assess the extent of this
potential efficiency loss.

The setup for the simulation study is similar to that of Hsu et
al. (2004) for investigating a semiparametric estimation of
marginal hazard function from case-control family study, with the
required modifications for the current prospective setting. For
each family we generated a common frailty value $W$ from the gamma
distribution with scale and shape parameters both equal
$\theta^{-1}$. We consider 300 families, each of size 2. A single
covariate from the standard normal distribution was incorporated.
Conditional on $W$, the survivor function is
 $$
 S(t|Z,W)=\exp\{-W \exp(\beta Z) (0.01t)^{4.6} \}
 $$
Thus, with $\beta=ln(2)$ or $ln(3)$ and a normal distribution for
the censoring, with mean 60 years and standard deviation of 15
years, the censoring level is approximately $85\%$ and $80\%$,
respectively. The censoring distribution was chosen in order to
generate appropriate mean age at onset and distribution, similar
to what is often observed for late onset diseases. With censoring
distributed according to $N(130,15^2)$ the respective censoring
levels are approximately $35\%$ and $30\%$. Table 1 summarizes the
results for the two estimation techniques, for $\beta^\circ=ln(2)$ or
$ln(3)$ and $\theta^\circ=2$. For our method, we compare the mean
estimated standard error based on our theoretical formula with the
empirical standard error, and provide the empirical coverage rate
of $95\%$ Wald-type confidence interval. For the EM-based method,
we report only the empirical standard error. In addition, the
empirical correlation between the EM-based estimators and our
estimators is presented. It is evident that both estimation
techniques perform very well in term of bias. Also, for our
method, good agreement was observed between the estimated and the
empirical standard error. The high values of the correlations
implies similarity between the two estimation techniques not only
on an average basis, but actually on a data set by data set basis.

\section{Example}
We now apply our method under the gamma frailty distribution to a
diabetic retinopathy data set. The Diabetic Retinopathy Study
(DRS) was begun in 1971 to study the effectiveness of laser
photocoagulation in delaying the onset of blindness in patients
with diabetic retinopathy. Patients with diabetic retinopathy and
visual acuity of 20/100 or better in both eyes were eligible for
the study. For each study subject, one eye was randomly selected
for treatment laser photocoagulation and the other eye was
observed without treatment. The outcome variable is time to
blindness of each eye. For illustrative purposes the following
analysis involves 197 high-risk patients as defined by DRS
criteria. Of the 394 measurements, 239 ($61\%$) are censored. The
regression coefficient estimate of the treatment effect was
$-0.890$ and $-0.910$ according to our proposed estimator and the
EM algorithm, respectively. The respective estimated standard
errors, $0.175$ and $0.174$, are based on 50 bootstrap samples.
The estimate of $\theta$ was $0.865$ with our approach, and
$0.857$ with the EM approach. The respective estimated standard
errors are $0.367$ and $0.365$.
As one can see, both method yield extremely similar
results. Both indicated that the treatment appeared effective in
delaying the time to blindness, and that the times to blindness
for both eyes are highly dependent. The hazard rate of one eye
becoming blind given the other eye is blind is almost twice
(1+$\theta$) as high as that given the other eye is not blind.

\section{Discussion}
We have presented a method for estimating the regression
coefficient vector and frailty parameter in a prospective frailty
survival model. The procedure is applicable to any frailty
distribution with finite moments. We have shown that our
estimators of the regression coefficients and frailty parameter
are consistent and asymptotically normally distributed, and given
an explicit consistent estimator for the variances of the
parameter estimates. For the popular gamma frailty model, we have
presented simulation results showing that our estimator is
essentially as efficient as the estimator based on the EM
algorithm. For our procedure, a consistent covariance estimator is
available which is much easier to compute than its counterpart for
the EM method as given by Parner (1998). Nonconjugate frailty
distributions can be handled by a simple univariate numerical
integration over the frailty distribution.

The estimation approach used here for estimating the cumulative
hazard function can be applied in some other important settings,
such as the case-control family study. Our approach avoids an
iterative procedure for the $\hat{\Lambda}_0$, enabling the
asymptotic properties of the estimator to be derived in a
relatively straightforward fashion. Shih and Chatterjee (2002)
proposed a semi-parametric quasi-partial-likelihood approach for
estimating the regression coefficients in survival data from a
case-control family study. Their cumulative hazard estimator
requires an iterative solution, and thus the properties of their
estimates could only be investigated so far by a simulation study.
If their method is modified by using our approach to estimating
$\Lambda_0(u)$, the proof presented in Appendix B can serve as a
basis for the asymptotic properties of the resulting procedure,
with appropriate modifications. The extension to this case will be
presented in a separate paper.

\section{Appendix: Asymptotic Theory}

\subsection{Assumptions and Background}

In deriving the asymptotic properties of $\hat{\bgamma}$ we make
the following assumptions:

\begin{enumerate}
\item
The random vectors
$(T^0_{i1},\ldots,T^0_{im_i},C_{i1},\ldots,C_{im_i},{\bf
Z}_{i1},\ldots,{\bf Z}_{im_i}, W_i)$, $i=1,\ldots,n$, are
independent and identically distributed.
\item
There is a finite maximum follow-up time $\tau>0$, with
$\E[\sum_{j=1}^{m_i}Y_{ij}(\tau)]=y^*>0$ for all $i$.
\item
\begin{enumerate}
\item
Conditional on ${\bf Z}_{ij}$ and $W_i$, the censoring is independent and
noninformative of $W_i$ and $(\bbeta,\Lambda_0)$.
\item
$W_i$ is independent of ${\bf Z}_{ij}$ and of $m_i$.\
\end{enumerate}
\item
The frailty random variable $W_i$ has finite moments up to order
$(m+2)$, where $m$ is a fixed upper bound on $m_{i}$.
\item
${\bf Z}_{ij}$ is bounded.
\item
The parameter $\bgamma$ lies in a compact subset $\cG$ of $\real^{p+1}$
containing an open neighborhood of $\bgamma^\circ$.
\item
There exist $b>0$ and $C>0$ such that
$$
\lim_{w \rightarrow 0} w^{-(b-1)} f(w) = C.
$$
\item
The baseline hazard function $\lambda_0^\circ(t)$ is bounded over $[0,\tau]$
by some constant $\lambda_{max}$.
\item
The function $f^\prime(w;\theta) = (d/d\theta) f(w;\theta)$ is absolutely
integrable.
\item
The censoring distribution has at most finitely many jumps on $[0,\tau]$.
\item
 The matrix $[(\partial/\partial \bgamma)
\bU(\bgamma,\hat{\Lambda}_0(\cdot,\bgamma))]|_{\bgams=\bgams^\circ}$
is invertible with probability going to 1 as $n \rightarrow
\infty$.
\end{enumerate}
The matrix
$(\partial/\partial\bgamma)\bU(\bgamma,\hat{\Lambda}_0(\cdot,\bgamma))$
is presented explicitly in Section 7.3, Step IV. From
(\ref{dls})-(\ref{dpp}), it is seen that a general proof of
invertibility is intractable, but given the data, one can easily
check that numerically the matrix is invertible.

\subsection{Technical Preliminaries}

Here we present some technical results that are needed for the
asymptotic theory. First note that,
by the boundedness of $\bbeta$ and ${\bf Z}_{ij}$,
there exists a constant $\nu>0$ such that
\begin{equation}
\nu^{-1} \leq \exp(\bbeta^T {\bf Z}_{ij}) \leq \nu.
\label{eb}
\end{equation}
Next, recall that $$ \psi_i( \bgamma,\Lambda,t) = \frac{\int
w^{N_i(t)+1} e^{-H_{i \cdot}(t)w} f(w) dw} {\int w^{N_i(t)}
e^{-H_{i \cdot}(t)w} f(w) dw}, $$ with $H_{i \cdot}(t)=H_{i
\cdot}(t,\bgamma,\Lambda) = \sum_{j=1}^{m_i} \Lambda(T_{ij} \wedge
t) \exp(\bbeta^T {\bf Z}_{ij})$ (here we define $H_{i \cdot}$ so
as to allow dependence on a general $\bgamma$ and $\Lambda$, which
will often not be explicitly indicated in the notation). Define
(for $0 \leq r \leq m$ and $h \geq 0$) $$ \psi^*(r,h) = \frac{\int
w^{r+1} e^{-hw} f(w) dw} {\int w^{r} e^{-hw} f(w) dw}. $$ Also
define $\psi_{min}^*(h) = \min_{0 \leq r \leq m} \psi^*(r,h)$ and
$\psi_{max}^*(h) = \max_{0 \leq r \leq m} \psi^*(r,h)$. Note that,
in the expression for $\psi^*(r,h)$, the numerator and denominator
are bounded from above by the assumption that $W$ has finite
$(m+2)$-th moment. In addition, the numerator and denominator are
by necessity strictly positive, for otherwise $W$ would have a
degenerate distribution concentrated at 0. Thus $\psi_{max}^*(h)$
is finite and $\psi_{min}^*(h)$ is strictly positive.

\vspace*{0.5em}
\noindent {\bf Lemma 1:} The function $\psi^*(r,h)$ is decreasing in $h$. In
consequence, we have, for all $\bgamma \in \cG$ and all $t$,
\begin{eqnarray}
\psi_i(\bgamma,\Lambda,t) & \leq & \psi_{max}^*(0), \label{huey} \\
\psi_i(\bgamma,\Lambda,t) & \geq & \psi_{min}^*(m \nu \Lambda(t)). \label{louie}
\end{eqnarray}
In addition, there exist $B>0$ and $\bar{h}>0$ such that,
for all $h \geq \bar{h}$,
\begin{equation}
\psi_{min}^*(h) \geq B h^{-1}.
\label{moe}
\end{equation}

\noindent {\bf Proof:}
We have
\begin{equation}
\frac{\partial}{\partial h} \psi^*(r,h)
= - \left[
\frac{\int w^{r+2} e^{-hw} f(w) dw}
{\int w^{r} e^{-hw} f(w) dw}
- \left (\frac{\int w^{r+1} e^{-hw} f(w) dw}
{\int w^{r} e^{-hw} f(w) dw} \right)^2 \right].
\label{psider}
\end{equation}
This quantity is negative for all $h$, which establishes that $\psi^*(r,h)$
is a decreasing function of $h$. Now, by definition, $\psi_i(\bgamma,\Lambda,t) =
\psi^*(N_i(t),H_{i \cdot}(t))$. We have $0 \leq H_{i \cdot}(t) \leq m \nu \Lambda(t)$.
The inequalities (\ref{huey}) and (\ref{louie}) follow immediately.

As for (\ref{moe}), using a change of variable and Assumption 7,
we find that $$ \lim_{h \rightarrow \infty} h \psi^*(r,h) =
\frac{\int_0^\infty v^{r+b} e^{-v} dv} {\int_0^\infty v^{r+b-1}
e^{-v} dv} = r + b. $$ Choosing $\bar{h}$ large enough so that the
above limit is obtained up to a factor of, say, 1.01, the result
follows.

\vspace*{0.5cm}
We define
$$
\bar{\Lambda} = 1.03 e^{m\sigma} \left( \frac{\bar{h}}{m \nu} \right),
$$
with $\sigma = 1.01 m \nu^2 /(By^*)$, where $\bar{h}$ and $B$ are as in Lemma 1.

\vspace*{0.5cm}
\noindent {\bf Lemma 2:}
With probability one, there exists $n^\prime$ such that, for all $t \in [0,\tau]$
and $\bgamma \in \cG$,
\begin{equation}
\hat{\Lambda}_0(t,\bgamma) \leq \bar{\Lambda} \quad \mbox{for } n \geq n^\prime.
\label{lamb}
\end{equation}

\noindent{\it Remark}: The point of this lemma is that $\hat{\Lambda}_0(t,\bgamma)$
is automatically bounded above, without any need to impose an upper bound artificially.

\noindent{\bf Proof:} To simplify the writing below, we will
suppress the argument $\bgamma$ in $\hat{\Lambda}_0(t,\bgamma)$.
Recall
$$
\Delta \hat{\Lambda}_0(\tau_k) =\frac
 {1}
 {\sum_{i=1}^{n}
 \psi_{i}(\bgamma,\hat{\Lambda}_0,\tau_{k-1})
 \sum_{j=1}^{m_i} Y_{ij}(\tau_k) \exp(\bbeta^T {\bf Z}_{ij})},
$$
where we now take $d_k=1$ since the survival time distribution is
assumed continuous. Using Lemma 1 and (\ref{eb}), we have
$$
\Delta \hat{\Lambda}_0(\tau_k) \leq
n^{-1} \nu
\psi_{min}^*(m \nu \hat{\Lambda}(\tau_{k-1}))
^{-1}
\left [ \frac{1}{n} \sum_{i=1}^n \sum_{j=1}^{m_i} Y_{ij}(\tau) \right]^{-1}.
$$
Now, since $\sum_{j=1}^{m_i}Y_{ij}(\tau)$ are iid random variables with
expectation $y^*$, by the strong law of large numbers we have
$$
\frac{1}{n} \sum_{i=1}^n \sum_{j=1}^{m_i} Y_{ij}(\tau)
\rightarrow y^*
$$
almost surely. Hence, with probability one, there exists $n^*$ such
that
\begin{equation}
\frac{1}{n} \sum_{i=1}^n \sum_{j=1}^{m_i} Y_{ij}(\tau)
\geq 0.999 y^*  \quad \mbox{for } n \geq n^*.
\label{yy}
\end{equation}
We thus have, for $n \geq n^*$,
\begin{equation}
\Delta \hat{\Lambda}_0(\tau_k) \leq
n^{-1} \left( \frac{1.01 \nu}{y^*} \right)
\psi_{min}^*(m \nu \hat{\Lambda}(\tau_{k-1}))^{-1}.
\label{del}
\end{equation}

Now, if $\hat{\Lambda}_0(t) \leq \bar{h}/(m\nu)$ for all $t$ then
we are done. Otherwise, there exists $k^\prime$ such that
$\hat{\Lambda}_0(\tau_k) \leq \bar{h}/(m\nu)$ for $k < k^\prime$
and $\hat{\Lambda}_0(\tau_k) \geq \bar{h}/(m\nu)$ for $k \geq
k^\prime$. Using the last inequality of Lemma 1, we obtain, for
$k>k^\prime$, $$ \Delta \hat{\Lambda}_0(\tau_k) \leq n^{-1} \sigma
\hat{\Lambda}_0(\tau_{k-1}), $$ or, in other words, $$
\hat{\Lambda}_0(\tau_k) \leq \left( 1 + \frac{\sigma}{n} \right)
\hat{\Lambda}_0(\tau_{k-1}). $$ Iterating the above inequality we
get $$ \hat{\Lambda}_0(\tau_{k^\prime+\ell}) \leq \left( 1 +
\frac{\sigma}{n} \right)^\ell \hat{\Lambda}_0(\tau_{k^\prime})
\leq \left( 1 + \frac{\sigma}{n} \right)^{mn}
\hat{\Lambda}_0(\tau_{k^\prime}) \leq 1.01 e^{m \sigma}
\hat{\Lambda}_0(\tau_{k^\prime}) $$ for $n$ large enough. But,
using (\ref{del}) and the fact that
$\hat{\Lambda}_0(\tau_{k^\prime-1}) \leq \bar{h}/(m\nu)$, we have
$$ \hat{\Lambda}_0(\tau_{k^\prime}) \leq \frac{\bar{h}}{m \nu} +
n^{-1} \left( \frac{1.01 \nu}{y^*} \right)
\psi_{min}^*(\bar{h})^{-1}, $$ which is less than $1.01
{\bar{h}}/({m \nu})$ for $n$ large enough. The desired conclusion
follows.

\vspace*{0.5cm}
\noindent {\bf Lemma 3:} $\sup_{s \in [0,\tau]}|\hat{\Lambda}_0(s,\bgamma^\circ)
-\hat{\Lambda}_0(s-,\bgamma^\circ)| \rightarrow 0 $ as $n \rightarrow \infty$.

\noindent{\bf Proof:}
Since $\hat{\Lambda}_0(s,\bgamma)-\hat{\Lambda}_0(s-,\bgamma)$
equals $\Delta\hat{\Lambda}_0(\tau_k)$ for $s=\tau_k$
and zero otherwise, it is enough to show
that $\sup_k \Delta \hat{\Lambda}_0(\tau_k,\bgamma^\circ)
\rightarrow 0$ as $n \rightarrow \infty$. But from Lemma 2 and
(\ref{del}) we have
$$
\Delta \hat{\Lambda}_0(\tau_k) \leq
n^{-1} \left( \frac{1.01 \nu}{y^*} \right)
\psi_{min}^*(m \nu \bar{\Lambda})^{-1}.
$$
for $n$ sufficiently large. The conclusion follows immediately.

\subsection{Consistency}

We now show the almost sure consistency of $\hat{\bbeta}$ and
$\hat{\Lambda}_0$. The argument is built on Claims~$\mbox{A-C}$ of
Section 3, which we prove below. Our argument follows
Zucker (2005, Appendix A.3).

\vspace*{0.5cm}
\noindent
{\bf Claim A:} $\hat{\Lambda}_0(t,\bgamma)$ converges a.s.\
to some function $\Lambda_0(t,\bgamma)$ uniformly in $t$ and $\bgamma$.

\noindent {\bf Proof:} In the proof below, whenever a functional norm is
written, the relevant uniform norm is intended.

Define $\Lambda_{max} = \max(\bar{\Lambda},\lambda_{max}\tau)$
and
$\psi^{**}(r,h) = \psi^*(r,h \wedge h_{max})$,
where $h_{max}$ = $m \nu \Lambda_{max}$. It is easy to see from (\ref{psider})
that $\psi^{**}(r,h)$ is Lipschitz continuous in $h$ (uniformly in~$r$). Recall
that $\psi_i(\bgamma,\Lambda,t) = \psi^*(N_i(t),H_{i \cdot}(t,\bgamma,\Lambda))$.
But Lemma 2 implies that $H_{i \cdot}(t,\bgamma,\hat{\Lambda}_0(\cdot,\bgamma)) \leq h_{max}$
for all $t \in [0,\tau]$ and $\bgamma \in \cG$. Hence we see that
$\psi_i(\bgamma,\hat{\Lambda}_0(\cdot,\bgamma),t) =
\psi^{**}(N_i(t),H_{i \cdot}(t,\bgamma,\hat{\Lambda}_0(\cdot,\bgamma)))$.

Now define, for a general function $\Lambda$,
$$
\Xi_n(t,\bgamma,\Lambda)=\int_0^t
\frac{n^{-1} \sum_{i=1}^{n}\sum_{j=1}^{m_i}dN_{ij}(s)}
{n^{-1} \sum_{i=1}^{n}\sum_{j=1}^{m_i}\psi^{**}(N_i(s-),H_{i \cdot}(s-,\bgamma,\Lambda))
Y_{ij}(s)\exp(\bbeta^T {\bf Z}_{ij})}
$$
and
$$
\Xi(t,\bgamma,\Lambda)=\int_0^t
 \frac{\E [\sum_{j=1}^{m_i}\psi^*(N_i(s-),H_{i \cdot}(s-,\bgamma^\circ,\Lambda_0^\circ))
Y_{ij}(s)\exp(\bbeta^{\circ T}{\bf Z}_{ij})]}
{\E [\sum_{j=1}^{m_i}\psi^{**}(N_i(s-),H_{i \cdot}(s-,\bgamma,\Lambda))Y_{ij}(s-)
\exp(\bbeta^T {\bf Z}_{ij})]} \lambda_0^\circ(s) ds.
$$
By definition, $\hat{\Lambda}_0(t,\bgamma)$ satisfies the equation
\begin{equation}
\hat{\Lambda}_0(t,\bgamma) = \Xi_n(t,\bgamma,\hat{\Lambda}_0(\cdot,\bgamma)).
\label{lest}
\end{equation}

Next, define
$$
q_{\bgams}(s,\Lambda)
= \frac{\E [\sum_{j=1}^{m_i}\psi^*(N_i(s-),H_{i \cdot}(s-,\bgamma^\circ,\Lambda_0^\circ))
Y_{ij}(s)\exp(\bbeta^{\circ T}{\bf Z}_{ij})]}
{\E [\sum_{j=1}^{m_i}\psi^{**}(N_i(s-),H_{i \cdot}(s-,\bgamma,\Lambda))Y_{ij}(s)
\exp(\bbeta^T {\bf Z}_{ij})]} \lambda_0^\circ(s).
$$
This function is uniformly bounded by
$B^* = [\psi_{max}^*(0)/\psi_{min}^*(h_{max})]\lambda_{max}$.
Moreover,
by the Lipschitz continuity of $\psi^{**}(r,h)$ with respect to $h$,
it satisfies the Lipschitz-like condition (for some constant $K$)
$$
|q_{\bgams}(s,\Lambda_1)-q_{\bgams}(s,\Lambda_2)| \leq
K \sup_{0 \leq u \leq s} |\Lambda_1(u)-\Lambda_2(u)|.
$$
Hence, by mimicking step by step the argument of Hartman (1973, Theorem 1.1), we find that
the equation $\Lambda(t) = \Xi(t,\bgamma,\Lambda)$ has a unique solution. We denote this
solution by $\Lambda_0(t,\bgamma)$. The claim then is that $\hat{\Lambda}_0(t,\bgamma)$
converges almost surely (uniformly in $t$ and $\bgamma$) to this function $\Lambda_0(t,\bgamma)$.
Though it may be possible to prove this claim directly, we shall use a convenient indirect argument.

Define $\tilde{\Lambda}_0^{(n)}(t,\bgamma)$ to be a modified version of
$\hat{\Lambda}_0(t,\bgamma)$ defined by linear interpolation between
the jumps, where we have added the superscript $n$ for emphasis.
Lemma 3 implies that, with probability one,
\begin{equation}
\sup_{t,\bgams} |\tilde{\Lambda}_0^{(n)}(t,\bgamma)
- \hat{\Lambda}_0(t,\bgamma)| \rightarrow 0, \label{Lone}
\end{equation}
and thus
\begin{equation}
\sup_{t,\bgams} |\Xi_n(t,\bgamma,\tilde{\Lambda}_0(t,\bgamma))-
\Xi_n(t,\bgamma,\hat{\Lambda}_0(t,\bgamma))| \rightarrow 0. \label{Ltwo}
\end{equation}
Lemma 2 shows that the family $\cL = \{ \tilde{\Lambda}_0^{(n)}(t,\bgamma)$,
$n \geq n^\prime \}$ is uniformly bounded. We will establish in a moment
that $\cL$ is also equicontinuous. It then follows, by the Arzela-Ascoli
theorem, that the closure of $\cL$ in $C([0,\tau] \times \cG)$ is compact.

The equicontinuity of $\cL$ is shown as follows.
Recall that $N_i(t)=\sum_{j=1}^{m_i} N_{ij}(t)$.
Write $\bar{N}(t) = n^{-1} \sum_{i=1}^n \sum_{j=1}^{m_i} N_{ij}(t)$.
We have $\bar{N}(t) \rightarrow \E[N_i(t)]$ as $n \rightarrow \infty$
uniformly in $t$ with probability one, with
$$
\E[N_i(t)] =
\int_0^t \E \left[ \sum_{j=1}^{m_i}\psi^*(N_i(s-),H_{i \cdot}(s-,\bgamma^\circ,\Lambda_0^\circ))
Y_{ij}(s)\exp(\bbeta^{\circ T}{\bf Z}_{ij}) \right] \lambda_0^\circ(s) ds.
$$
In view of this and (\ref{yy}) there exists a probability-one set of realizations
$\Omega^*$ on which the following holds: for any given $\epsilon>0$, we can find
$n^{\prime \prime}(\epsilon)$ such that $\sup_t |\bar{N}(t)-E[N_i(t)]| \leq \epsilon/(4B^\circ)$
for all $n \geq n^{\prime \prime}(\epsilon)$, where $B^\circ = 1.01 \nu/
[\psi_{min}^*(h_{max})y^*]$.
In consequence, for all $t$ and $u$ with $u<t$, we find that
$$
 \hat{\Lambda}_0(t,\bgamma) - \hat{\Lambda}_0(u,\bgamma)
 = \int_u^t
\frac{n^{-1} \sum_{i=1}^{n}\sum_{j=1}^{m_i}dN_{ij}(s)}
{n^{-1} \sum_{i=1}^{n}\sum_{j=1}^{m_i}\psi^{**}(N_i(s-),H_{i \cdot}(s-,\bgamma,\Lambda))
Y_{ij}(s)\exp(\bbeta^T {\bf Z}_{ij})}
$$
satisfies
\begin{equation}
\hat{\Lambda}_0(t,\bgamma) - \hat{\Lambda}_0(u,\bgamma)
\leq B^* (t-u) + \frac{\epsilon}{2} \quad \mbox{for all } n \geq n^{\prime \prime}(\epsilon).
\label{equi}
\end{equation}
Moreover, it is easy to see that $\hat{\Lambda}_0(t,\bgamma)$ is
Lipschitz continuous in $\bgamma$ with Lipschitz constant
$C^*$, say, that is independent of $t$.

These two results imply that $\cL$ is equicontinuous. This is seen as follows.
For given $\epsilon$, we need to find $\delta_1^*$ and $\delta_2^*$ such that
$|\tilde{\Lambda}_0^{(n)}(t,\bgamma) -\tilde{\Lambda}_0^{(n)}(u,\bgamma)| \leq \epsilon$ whenever $|t-u|
\leq \delta_1^*$ and $|\tilde{\Lambda}_0^{(n)}(t,\bgamma) -\tilde{\Lambda}_0^{(n)}(t,\bgamma^\prime)| \leq
\epsilon$ whenever $\|\bgamma-\bgamma^\prime \| \leq \delta_2^*$. The latter is easily obtained using the
Lipschitz continuity of $\hat{\Lambda}_0(t,\bgamma)$ with respect to $\bgamma$. As for the former,
for $n \geq n^{\prime \prime}(\epsilon)$ this can be accomplished using (\ref{equi}), while for $n$ in
the finite set $n^\prime \leq n < n^{\prime \prime}(\epsilon)$ this can be accomplished using the
fact that the function $\tilde{\Lambda}_0^{(n)}(t,\bgamma)$ is uniformly continuous on $[0,\tau]$ for
every given $n$.

We have thus shown that $\cL$ is (almost surely)
a relatively compact set in the space $C([0,\tau] \times \cG)$.

Next, define
\begin{eqnarray*}
A(\bgamma,\Lambda,s) & = &
\frac{1}{n} \sum_{i=1}^{n}\sum_{j=1}^{m_i}\psi^{**}(N_i(s-),H_{i \cdot}(s-,\bgamma,\Lambda))
Y_{ij}(s)\exp(\bbeta^T {\bf Z}_{ij}), \\
a(\bgamma,\Lambda,s) & = &
\E \left[\sum_{j=1}^{m_i}\psi^{**}(N_i(s-),H_{i \cdot}(s-,\bgamma,\Lambda))
Y_{ij}(s)\exp(\bbeta^T {\bf Z}_{ij}) \right].
\end{eqnarray*}
We show below that, with probability one,
\begin{equation}
\sup_{s,\bgams} |A(\bgamma,\tilde{\Lambda}^{(n)},s) - a(\bgamma,\tilde{\Lambda}^{(n)},s)| \rightarrow 0.
\label{acn}
\end{equation}
Given this and the a.s. uniform convergence of
$\bar{N}(t)$ to $\E[N_i(t)]$, we can infer that
\begin{equation}
\sup_{t,\bgams}
|\Xi_n(t,\bgamma,\tilde{\Lambda}_0^{(n)}(t,\bgamma))
- \Xi(t,\bgamma,\tilde{\Lambda}_0^{(n)}(t,\bgamma))| \rightarrow 0.
\label{supX}
\end{equation}
The result (\ref{supX}) is easily obtained by adapting the
argument of Aalen (1976, Lemma~6.1), making use of the
equicontinuity of $\cL$. It is here that we make use of
Assumption~10, for the adaptation of Aalen's argument requires
$a(\bgamma,\Lambda,s)$ to be piecewise continuous with finite left
and right limits at each point of discontinuity.

From (\ref{lest}), (\ref{Lone}), (\ref{Ltwo}), and (\ref{supX}) it
follows that any limit point of $\{
\tilde{\Lambda}_0^{(n)}(t,\bgamma) \}$ must satisfy the equation
$\Lambda=\Xi(t,\bgamma,\Lambda)$. Since $\Lambda_0(t,\bgamma)$ is
the unique solution of this equation, it is the unique limit point
of $\{ \tilde{\Lambda}_0^{(n)}(t,\bgamma) \}$. Thus $\{
\tilde{\Lambda}_0^{(n)}(t,\bgamma) \}$ is a sequence in a compact
set with unique limit point $\Lambda_0(t,\bgamma)$. Hence
$\tilde{\Lambda}_0^{(n)}(t,\bgamma)$ converges a.s.\ uniformly in
$t$ and $\bgamma$ to $\Lambda_0(t,\bgamma)$. In view of
(\ref{Lone}), the same holds of $\hat{\Lambda}_0(t,\bgamma)$,
which is the desired result.

To complete the proof, we must establish (\ref{acn}). This involves several steps.
First, it is easy to see that there exists a constant $\kappa$ (independent of
$\bgamma$ and $s$) such that
\begin{eqnarray}
\sup_{s,\bgams} | A(\bgamma,\Lambda_1,s) - A(\bgamma,\Lambda_2,s) |
& \leq & \kappa \| \Lambda_1 - \Lambda_2 \|, \label{Alip} \\
\sup_{s,\bgams} | a(\bgamma,\Lambda_1,s) - a(\bgamma,\Lambda_2,s) |
& \leq & \kappa \| \Lambda_1 - \Lambda_2 \|. \label{alip}
\end{eqnarray}
Next, for any fixed continuous $\Lambda$, the functional strong law of large numbers of
Andersen \& Gill (1982, Appendix III) implies that, with probability one,
\begin{equation}
\sup_{s,\bgams} |A(\bgamma,\Lambda,s) - a(\bgamma,\Lambda,s)| \rightarrow 0.
\label{ac}
\end{equation}

Now, given $\eps>0$, define the sets $\{ t_j^\peps \}$, $\{ \bgamma_k^\peps \}$,
and $\{ \Lambda_l^\peps \}$
to be finite partition grids of $[0,\tau]$, $\cG$, and
$[0,\Lambda_{max}]$, respectively, with distance of no more than $\eps$ between grid
points. Define $\cL_\eps^*$ to be the set of functions of $t$ and $\bgamma$ defined by
linear interpolation through vertices of the form $(t_j^\peps,\bgamma_k^\peps,\Lambda_l^\peps)$.

Obviously $\cL_{\eps}^*$ is a finite set. Hence, in view of (\ref{ac}), there exists
a probability-one set of realizations $\Omega_\eps$ for which
\begin{equation}
\sup_{s \in [0,\tau], \bgams \in \cG, \Lambda \in \cL_\eps^*}
|A(\bgamma,\Lambda,s) - a(\bgamma,\Lambda,s)| \rightarrow 0.
\label{acle}
\end{equation}
Define
$$
\Omega^{**} =
\bigcap_{\ell=1}^\infty \Omega_{1/\ell}
$$
and $\Omega_0 = \Omega^* \cap \Omega^{**}$, with $\Omega^*$ as defined earlier.
Clearly $\Pr(\Omega_0)=1$. From now on, we restrict attention to $\Omega_0$.

Now let $\eps>0$ be given. Choose $\ell > \eps^{-1}$. In view of (\ref{equi}) and
(\ref{acle}), we can find for any $\omega \in \Omega_0$ a suitable positive integer
$\bar{n}(\eps,\omega)$ such that, whenever $n \geq \bar{n}(\eps,\omega)$,
\begin{equation}
|\tilde{\Lambda}_0^{(n)}(t,\bgamma) - \tilde{\Lambda}_0^{(n)}(u,\bgamma)|
\leq B^* (t-u) + \frac{\epsilon}{2} \quad \forall t,u,
\label{thingone}
\end{equation}
\begin{equation}
\sup_{s \in [0,\tau], \bgams \in \cG, \Lambda \in \cL_{1/\ell}^*}
|A(\bgamma,\Lambda,s) - a(\bgamma,\Lambda,s)| \leq \eps.
\label{thingtwo}
\end{equation}

Next, let $\bar{\Lambda}_0^{(n)}$ denote the function defined by linear interpolation
through $(t_j^\peps,\bgamma_k^\peps,\bar{\Lambda}_{jk}^\peps)$, where $\bar{\Lambda}_{jk}^\peps$
is the element of $\{ \Lambda_l^\peps \}$ that is closest to $\tilde{\Lambda}_0^{(n)}
(t_j^\peps,\bgamma_k^\peps)$. It is clear that
$$
|\bar{\Lambda}_0^{(n)}(t_j^\peps,\bgamma_k^\peps)
- \tilde{\Lambda}_0^{(n)}(t_j^\peps,\bgamma_k^\peps)| \leq \eps \quad \forall j,k.
$$
Using (\ref{thingone}) and the Lipschitz continuity of
$\tilde{\Lambda}_0^{(n)}(t,\bgamma)$ with respect to $\bgamma$
(which follows from the corresponding property of $\hat{\Lambda}_0(t,\bgamma))$,
we thus obtain
$$
\sup_{t,\bgams} |\bar{\Lambda}_0^{(n)}(t,\bgamma) - \tilde{\Lambda}_0^{(n)}(t,\bgamma)|
\leq B^{**} \eps
$$
for a suitable fixed constant $B^{**}$ (depending on $B^*$ and $C^*$).
Combining this with (\ref{thingtwo}) and (\ref{alip}), we obtain
$$
\sup_{s,\bgams} |A(\bgamma,\tilde{\Lambda}^{(n)},s) - a(\bgamma,\tilde{\Lambda}^{(n)},s)|
\leq (2\kappa B^{**} + 1) \eps \quad \mbox{for all } n \geq \bar{n}(\eps,\omega).
$$
Since $\eps$ was arbitrary, the desired conclusion (\ref{acn}) follows, and the proof is
thus complete.
\\

\noindent
{\it Remark}: Note that $\Lambda_0(\cdot,\bgamma^\circ)=\Lambda_0^\circ(\cdot)$
since $\Lambda_0^\circ$ trivially solves the equation $\Lambda=\Xi(t,\bgamma^\circ,\Lambda)$.
\\

\noindent
{\bf Claim B:} With probability one, $\bU(\bgamma,\hat{\Lambda}_0(\cdot,\bgamma))$
converges to $\bu(\bgamma,\Lambda_0(\cdot,\bgamma))
= \E[ \bU(\bgamma,\Lambda_0(\cdot,\bgamma)) ]$ uniformly in $\bgamma$ over $\cG$.

\noindent {\bf Proof:}
Since $\bU(\bgamma,\Lambda_0(\cdot,\bgamma))$ is the mean of iid terms,
the functional strong law of numbers of Andersen \& Gill (1982, Appendix III) implies that
$\bU(\bgamma,\Lambda_0(\cdot,\bgamma))$ converges uniformly in $\bgamma$ almost surely
to $\bu(\bgamma,\Lambda_0(\cdot,\bgamma))$. It remains only to show that
\begin{equation}
\sup_{\bgams} |\bU(\bgamma,\hat{\Lambda}_0(\cdot,\bgamma)) - \bU(\bgamma,\Lambda_0(\cdot,\bgamma))|
\rightarrow 0
\label{fuzz}
\end{equation}
almost surely. Now it may be seen easily from the structure of $\bU(\bgamma,\Lambda)$ that there
exists some constant $C^\circ$ (independent of $\bgamma$) such that
$$
|\bU(\bgamma,\Lambda_1)-\bU(\bgamma,\Lambda_2)| \leq C^\circ \| \Lambda_1 - \Lambda_2 \|.
$$
Given this result along with the result of Claim A, the result
(\ref{fuzz}) follows immediately.
\\

\noindent {\bf Claim C:} There exists a unique consistent root to
${\bf U}(\hat{\bgamma},\hat{\Lambda}_0(\cdot,\hat{\bgamma}))=0$.

\noindent {\bf Proof:}
We apply Foutz's (1977) theorem on consistency of maximum likelihood type
estimators. The following conditions must be verified:
\begin{description}
\item[F1.]
$\partial \bU(\bgamma,\hat{\Lambda}_0(\cdot,\bgamma))/\partial
\bgamma$ exists and is continuous in an open neighborhood about
$\bgamma^\circ$.
\item[F2.]
The convergence of $\partial
\bU(\bgamma,\hat{\Lambda}_0(\cdot,\bgamma)) /\partial \bgamma$ to
its limit is uniform in open neighborhood of $\bgamma^\circ$.
\item[F3.]
$\bU(\bgamma^\circ,\hat{\Lambda}_0(\cdot,\bgamma^\circ)) \rightarrow 0$ as
$n \rightarrow \infty$.
\item[F4.]
The matrix $-[\partial
\bU(\bgamma,\hat{\Lambda}_0(\cdot,\bgamma))/\partial
\bgamma]|_{\bgams=\bgams^\circ}$ is invertible with probability
going to 1 as $n \rightarrow \infty$. (In Foutz's paper, the
matrix in question is symmetric, and so he stated the condition in
terms of positive definiteness. But it is clear from his proof,
which is based on the inverse function theorem, that the basic
condition needed is invertibility.)
\end{description}
It is easily seen that Condition F1 holds. Given Assumptions 2, 4,
and 5, Condition F2 follows from the previously-cited functional
law of large numbers. As for Condition F3, in Claim B we showed
that $\bU(\bgamma,\Lambda_0(\cdot,\bgamma))$ converges a.s.\
uniformly to $\bu(\bgamma,\Lambda_0(\cdot,\bgamma)) = \E
[\bU(\bgamma,\Lambda_0(\cdot,\bgamma))]$. We noted already that
$\Lambda_0(\cdot,\bgamma^\circ)=\Lambda_0(\cdot)$. Thus all we
need is to show that $\E [\bU(\bgamma^\circ,\Lambda_0)]=~{\bf 0}$.
Since $\bU$ is a score function derived from a classical iid
likelihood, this result follows from classical likelihood theory.
Condition F4 has been assumed in Assumption 11; we noted
previously that, given the data, it can be checked numerically.
With Conditions F1-F4 thus verified, it follows from Foutz's
theorem that $\hat{\bgamma} \rightarrow \bgamma^\circ$ as $n
\rightarrow \infty$ with probability one.

\subsection{Asymptotic Normality}
To show that $\hat{\bgamma}$ is asymptotically normally
distributed, we write
\begin{eqnarray}
\lefteqn{{\bf 0}}&=&\bU(\hat{\bgamma},\hat{\Lambda}_0(\cdot,\hat{\bgamma}))\nonumber
\\ &=& \bU(\bgamma^\circ,\Lambda_0^\circ) +
[\bU(\bgamma^\circ,\hat{\Lambda}_0(\cdot,\bgamma^\circ))-\bU(\bgamma^\circ,\Lambda_0^\circ)]
\nonumber\\ & & + \,
[\bU(\hat{\bgamma},\hat{\Lambda}_0(\cdot,\hat{\bgamma}))-
\bU(\bgamma^\circ,\hat{\Lambda}_0(\cdot,\bgamma^\circ))] \nonumber
\end{eqnarray}
In the following we consider each of the above terms of the
right-hand side of the equation.


\bsh
\underline{Step I}
\esh

We can write $\bU(\bgamma^\circ, \Lambda_0^\circ)$ as
$$
\bU(\bgamma^\circ, \Lambda_0^\circ)
= \frac{1}{n} \sum_{i=1}^n \bxi_i,
$$
where $\bxi_i$ is a $(p+1)$-vector with $r$-th element, $r=1,\ldots,p$,
given by
$$
\xi_{ir}  =  \sum_{j=1}^{m_i} \delta_{ij} Z_{ijr} -
 \frac{\left[\sum_{j=1}^{m_i}{H}_{ij}(\tau)
 Z_{ijr}\right] \int w^{N_{i.}(\tau)+1} \exp\{-w \{H_{i.}(\tau)\}f(w;{\theta})dw}
 {\int w^{N_{i.}(\tau)} \exp\{-w {H}_{i.}(\tau)\}f(w;{\theta})dw}
$$
and $(p+1)$-th element given by
$$
\xi_{i(p+1)} =  \frac{\int w^{N_{i.}(\tau)}
\exp\{-w{H}_{i.}(\tau)\} f'(w;{\theta})dw}
 {\int w^{N_{i.}(\tau)} \exp\{-w {H}_{i.}(\tau)\}f(w;{\theta})dw}.
$$
Thus $\bU(\bgamma^\circ, \Lambda_0^\circ)$ is the
mean of the iid mean-zero random vectors $\bxi_i$. It hence follows immediately
from the classical central limit theorem that $n^{\half} \bU(\bgamma^\circ, \Lambda_0^\circ)$
is asymptotically mean-zero multivariate normal. To estimate the covariance
matrix, let $\bxi_i^*$ be the counterpart of $\bxi_i$ with estimates of $\bgamma$
and $\Lambda_0$ substituted for the true values. Then an empirical
estimator of the covariance matrix is given by
 $$
 \hat{{\bf V}}(\hat{\bgamma}) = \frac{1}{n} \sum_{i=1}^{n}
 \bxi_i^{*}\bxi_i^{*T}.
 $$
This is a consistent estimator of the covariance matrix
since $\hat{\Lambda}_0(t,\bgamma)$ converges to $\Lambda_0(t,\bgamma)$ a.s.\
uniformly in $t$ and $\bgamma$  (Claim A), and $\hat{\bgamma}$ is a
consistent estimator of $\bgamma^\circ$ (Claim C).


\bsh
\underline{Step II}
\esh

Let $\hat{U}_r=U_r(\bgamma^\circ,\hat{\Lambda}_0)$, $r=1,\ldots,p$,
and $\hat{U}_{p+1}=U_{p+1}(\bgamma^\circ,\hat{\Lambda}_0)$ (in this
segment of the proof, when we write $(\bgamma^\circ,\hat{\Lambda}_0)$
the intent is to signify $(\bgamma^\circ,\hat{\Lambda}_0(\cdot,\bgamma^\circ))$.
First order Taylor expansion of $\hat{U}_r$ about $\Lambda_0^\circ$,
$r=1,\ldots,p+1$, gives
$$
\hspace*{-6cm}
n^{1/2}\{U_r(\bgamma^\circ,\hat{\Lambda}_0)-U_r(\bgamma^\circ,{\Lambda}_0^\circ)\}
$$
\begin{equation}\label{eq:tayloru}
\hspace*{1cm} = n^{-1/2} \sum_{i=1}^{n} \sum_{j=1}^{m_i}
Q_{ijr}(\bgamma^\circ,\Lambda^\circ,T_{ij})
\{\hat{\Lambda}_0(T_{ij},\bgamma^\circ)-\Lambda_0^\circ(T_{ij})\} +
o_p(1),
\end{equation}
where
\begin{eqnarray*}
Q_{ijr}(\bgamma^\circ,\Lambda^\circ,T_{ij})&=& -\left\{
 \frac{\phi_{2i}(\bgamma^\circ,\Lambda_0^\circ,\tau)}{\phi_{1i}(\bgamma^\circ,\Lambda_0^\circ,\tau)}
 R_{ij}^* Z_{ijr}
 -
 \frac{\phi_{3i}(\bgamma^\circ,\Lambda_0^\circ,\tau)}{\phi_{1i}(\bgamma^\circ,\Lambda_0^\circ,\tau)}
 R_{ij}^*\sum_{j=1}^{m_i} H_{ij}(T_{ij}) Z_{ijr}
 \right. \nonumber \\
& & \left.
\hspace*{1cm} + \, \frac{\phi_{2i}^2(\bgamma^\circ,\Lambda_0^\circ,\tau)}
 {\phi_{1i}^2(\bgamma^\circ,\Lambda_0^\circ,\tau)}
 R_{ij}^* \sum_{j=1}^{m_i} H_{ij}(T_{ij})Z_{ijr}  \right\}
\end{eqnarray*}
for $r=1,\ldots,p$, and
\begin{eqnarray*}
Q_{ij(p+1)}(\bgamma^\circ,\Lambda^\circ,T_{ij}) =  R_{ij}^* \left\{
\frac{\phi_{2i}(\bgamma^\circ,\Lambda_0^\circ,\tau)
 \phi_{1i}^{(\theta)}(\bgamma^\circ,\Lambda_0^\circ,\tau)}
 {\phi_{1i}^2(\bgamma^\circ,\Lambda_0^\circ,\tau)}
 -\frac{ \phi_{2i}^{(\theta)}(\bgamma^\circ,\Lambda_0^\circ,\tau) }
 {\phi_{1i}(\bgamma^\circ,\Lambda_0^\circ,\tau)}
\right\},
\end{eqnarray*}
with $R_{ij}^*=\exp(\bbeta^T {\bf Z}_{ij})$ and
 \begin{eqnarray*}
 \phi_{ki}^{(\theta)}(\bgamma,\Lambda_0,t)
 =
 \int w^{N_{i.}(t)+(k-1)} \exp\{-w H_{i.}(t)\} f'(w) dw,
\quad k=1,2.
 \end{eqnarray*}
The validity of the approximation (\ref{eq:tayloru}) can be seen by an argument
similar to that used in connection with (\ref{eq:yapp}) below.

Based on the intensity process (\ref{eq:inten}), the process
 \begin{eqnarray*}
 M_{ij}(t)=N_{ij}(t)-\int_0^t
 \lambda_0(u) \exp(\bbeta^{\circ T} {\bf Z}_{ij}) Y_{ij}(u) \psi_{i}(\bgamma^\circ,\Lambda_0^\circ,u-) du
 \end{eqnarray*}
is a mean zero martingale with respect to the filtration
${\mathcal F}_t$. Also, by Lemma 3, we have that
$\sup_{s \in [0,\tau]} |\hat{\Lambda}_0(s,\bgamma^\circ)-\hat{\Lambda}_0(s-,\bgamma^\circ)|$
converges to zero. Thus, replacing $s-$ by $s$ we obtain the following
approximation, uniformly over $t \in [0,\tau]$:
 \begin{eqnarray}\label{eq:lamapp}
\hat{\Lambda}_0(t,\bgamma^\circ)-\Lambda_0^\circ(t) & \approx & \frac{1}{n}
\int_0^t \{\mathcal{Y}(s,\Lambda_0^\circ)\}^{-1} \sum_{i=1}^{n}
 \sum_{j=1}^{m_i} dM_{ij}(s) \nonumber\\
 &+& \frac{1}{n} \int_0^t \left[\{\mathcal{Y}(s,\hat{\Lambda}_0)\}^{-1}
 - \{\mathcal{Y}(s,\Lambda_0^\circ)\}^{-1}\right]\sum_{i=1}^{n}
 \sum_{j=1}^{m_i} dN_{ij}(s),
 \end{eqnarray}
where
 $$
 \mathcal{Y}(s,\Lambda) = \frac{1}{n}  \sum_{i=1}^{n}
 \psi_i(\bgamma^\circ,\Lambda,s)\sum_{j=1}^{m_i} Y_{ij}(s)
 \exp(\bbeta^{\circ T} {\bf Z}_{ij}).
 $$

Now let
 $$
 {\mathcal W}(s,r)=\{{\mathcal
 Y}(s,\Lambda_0^\circ+r\Delta)\}^{-1}
 $$
with $\Delta=\hat{\Lambda}_0-\Lambda_0^\circ$. Define
$\dot{\mathcal W}$ and $\ddot{\mathcal W}$ as the first and second
derivative of ${\mathcal W}$ with respect to $r$, respectively.
Then, by a first order Taylor expansion of ${\mathcal W}(s,r)$
around $r=0$ evaluated at $r=1$ with Lagrange remainder
(Abramowitz and Stegun, 1972, p. 880)  we get (after computing the
necessary derivatives)
\begin{eqnarray}
\lefteqn{\{{\mathcal Y}(s,\hat{\Lambda}_0)\}^{-1} - \{{\mathcal
Y}(s,{\Lambda}_0^\circ)\}^{-1}
 =\dot{\mathcal
W}(s,0)+\frac{1}{2}\ddot{\mathcal W}(s,\tilde{r}(s))}
\nonumber \\
&& \hspace*{-2cm} =  -\frac{1}{n} \sum_{i=1}^{n}\sum_{j=1}^{m_i} \left[
 \frac{R_{i.}(s)\eta_{1i}(0,s)}{\{{\mathcal
Y}(s,{\Lambda}_0^\circ)\}^2}-\frac{1}{2}h_i(\tilde{r}(s),s)
  \right]
  \exp(\bbeta^T {\bf Z}_{ij})
  \{ \hat{\Lambda}_0(T_{ij}\wedge s) - {\Lambda}_0^\circ(T_{ij}\wedge s)\},
  \label{eq:yapp}
\end{eqnarray}
where $R_{ij}(u)=\exp(\bbeta^T {\bf Z}_{ij}) Y_{ij}(u)$,
$R_{i.}(u)=\sum_{j=1}^{m_i} R_{ij}(u)$, $\tilde{r}(s) \in [0,1]$,
\begin{eqnarray*}
\eta_{1i}(r,s)=\frac{\phi_{3i}(\bgamma^\circ,{\Lambda}_0^\circ+r\Delta,s)}
{\phi_{1i}(\bgamma^\circ,{\Lambda}_0^\circ+r\Delta,s)}-\left\{
\frac{\phi_{2i}(\bgamma^\circ,{\Lambda}_0^\circ+r\Delta,s)}{\phi_{1i}(\bgamma^\circ,{\Lambda}_0^\circ+r\Delta,s)}
\right\}^2,
\end{eqnarray*}
and $h_i(r,s)$ is as defined \S 7.5 below. We show there that $h_i(r,s)$
is $o(1)$ uniformly in $r$ and $s$.

Let $\eta_{1i}(s)=\eta_{1i}(0,s)$. Plugging (\ref{eq:yapp}) into
(\ref{eq:lamapp})
 we get
 \begin{eqnarray*}
 \lefteqn{\hat{\Lambda}_0(t,\bgamma^\circ)-\Lambda_0^\circ(t)  \approx
 n^{-1}\int_0^t \{\mathcal{Y}(s,\Lambda_0^\circ)\}^{-1} \sum_{i=1}^{n}
 \sum_{j=1}^{m_i} dM_{ij}(s)} \\
 && -n^{-2}\int_0^t
 \sum_{k=1}^{n}\sum_{l=1}^{m_k}
 \frac{I(T_{kl}>s) R_{k.}(s) \eta_{1k}(s)}
 {\{\mathcal{Y}(s,\Lambda_0^\circ)\}^{2}}
 \exp(\bbeta^T {\bf Z}_{kl})
 \{\hat{\Lambda}_0(s)-\Lambda_0^\circ(s)\}\sum_{i=1}^{n}
 \sum_{j=1}^{m_i} dN_{ij}(s)  \\
 && - n^{-2}\int_0^t
\sum_{k=1}^{n}\sum_{l=1}^{m_k}\frac{I(T_{kl} \leq s) R_{k.}(s)
 \eta_{1k}(s)}
 {\{\mathcal{Y}(s,\Lambda_0^\circ)\}^{2}}
 \exp(\bbeta^T {\bf Z}_{kl})
 \{\hat{\Lambda}_0(T_{kl})-\Lambda_0^\circ(T_{kl})\}\sum_{i=1}^{n}
 \sum_{j=1}^{m_i} dN_{ij}(s) \\
 && + \, n^{-2}\int_0^t
 \sum_{k=1}^{n}\sum_{l=1}^{m_k} \frac{1}{2} h_{k}(\tilde{r}(s),s) \exp(\bbeta^T {\bf Z}_{kl})
 \{\hat{\Lambda}_0(T_{kl})-\Lambda_0^\circ(T_{kl})\}\sum_{i=1}^{n}
 \sum_{j=1}^{m_i} dN_{ij}(s).
\end{eqnarray*}
The third term of the above equation can be written, by interchanging the
order of integration, as
$$
n^{-2}\sum_{k=1}^{n}\sum_{l=1}^{m_k}\sum_{i=1}^{n}\sum_{j=1}^{m_i}
\int_0^t \frac{R_{k.}(s) \eta_{1k}(s)}
{\{\mathcal{Y}(s,\Lambda_0^\circ)\}^{2}} \exp(\bbeta^T {\bf Z}_{kl})
 \left[ \int_0^s
\{\hat{\Lambda}_0(u)-\Lambda_0^\circ(u) \} d\tilde{N}_{kl}(u)\}
\right] dN_{ij}(s)
$$
$$
 = \int_0^{t}
  \{\hat{\Lambda}_0(s)-\Lambda_0^\circ(s)\}
\sum_{i=1}^{n}\sum_{j=1}^{m_i} \Omega_{ij}(s,t) d\tilde{N}_{ij}(s),
$$
where
$\tilde{N}_{ij}(t)=I(T_{ij} \leq t)$ and
 $$
 \Omega_{ij}(s,t)=n^{-2}\int_s^t \{\mathcal{Y}(u,\Lambda_0^\circ)\}^{-2}
 R_{i.}(u) \eta_{1i}(u)\exp(\bbeta^T {\bf Z}_{ij}) \sum_{k=1}^n \sum_{l=1}^{m_k}
 dN_{kl}(u).
 $$
Hence we get
 \begin{eqnarray}\label{eq:lamapp2}
\hat{\Lambda}_0(t,\bgamma^\circ)-\Lambda_0^\circ(t)  &\approx&
 n^{-1}\int_0^t \{\mathcal{Y}(s,\Lambda_0^\circ)\}^{-1} \sum_{i=1}^{n}
 \sum_{j=1}^{m_i} dM_{ij}(s) \nonumber\\
 &-&  \int_0^t
 \{\hat{\Lambda}_0(s,\bgamma^\circ)-\Lambda_0^\circ(s)\}\sum_{i=1}^{n}
 \sum_{j=1}^{m_i}\{\delta_{ij} \Upsilon(s)+\Omega_{ij}(s,t) +o(n^{-1})\}
 d\tilde{N}_{ij}(s) \nonumber
 \end{eqnarray}
where
 $$
\Upsilon(s) =n^{-2}\{\mathcal{Y}(s,\Lambda_0^\circ)\}^{-2}
\sum_{k=1}^{n} \sum_{l=1}^{m_k} I(T_{kl}>s) R_{k.}(s)
\eta_{1k}(s)\exp(\bbeta^T {\bf Z}_{kl}).
 $$
The $o(n^{-1})$ is uniform in $t$ (see \S 7.5) and will be
dominated by $\Omega$ and $\Upsilon$, which are of order $n^{-1}$.
Hence the $o(n^{-1})$ term can be ignored.

Given the all the above, an argument similar to that
of Yang \& Prentice (1999) and Zucker (2005) yields
following martingale representation
\begin{eqnarray}\label{eq:martin}
\hat{\Lambda}_0(t,\bgamma^\circ)-\Lambda_0^\circ(t)
 &\approx&
 \frac{1}{n \hat{p}(t)} \int_0^t
 \frac{\hat{p}(s-) \sum_{i=1}^{n}\sum_{j=1}^{m_i} dM_{ij}(s)}
 {\mathcal{Y}(s,\Lambda_0^\circ)},
\end{eqnarray}
where
\begin{eqnarray*}
\hat{p}(t)=\prod_{s \leq t} \left[
 1+\sum_{i=1}^{n} \sum_{j=1}^{m_i}
 \{\delta_{ij} \Upsilon(s)+\Omega_{ij}(s,t) \} d\tilde{N}_{ij}(s)
 \right].
\end{eqnarray*}
Based on (\ref{eq:tayloru}), we can write
 \begin{eqnarray*}
U_r(\bgamma^\circ,\hat{\Lambda}_0)-U_r(\bgamma^\circ,{\Lambda}_0^\circ)
\approx n^{-1} \sum_{i=1}^{n} \sum_{j=1}^{m_i}
 \int_0^{\tau} Q_{ijr}(\bgamma^\circ,\Lambda_0^\circ,s)
\{\hat{\Lambda}_0(s,\bgamma^\circ)-\Lambda_0^\circ(s)\} d
\tilde{N}_{ij}(s).
\end{eqnarray*}
Plugging the martingale representation (\ref{eq:martin}) into the
above equation and interchanging the order of integration gives
\begin{eqnarray}
U_r(\bgamma^\circ,\hat{\Lambda}_0)-U_r(\bgamma^\circ,{\Lambda}_0^\circ)
& & \nonumber \\
& \hspace*{-6cm} \approx & \hspace*{-3cm}  n^{-2} \sum_{i=1}^{n} \sum_{j=1}^{m_i}  \int_0^\tau
 \frac{Q_{ijr}(\bgamma^\circ,\Lambda_0^\circ,t)}{\hat{p}(t)} \int_0^t
 \frac{\hat{p}(s-) \sum_{k=1}^{n}\sum_{l=1}^{m_k} dM_{kl}(s)}
 {\mathcal{Y}(s,\Lambda_0^\circ)}
 d\tilde{N}_{ij}(t)
\nonumber \\
& \hspace*{-6cm}  =  & \hspace*{-3cm}
n^{-1}
\int_0^{\tau} \pi_r(s,\bgamma^\circ,\Lambda_0^\circ) \frac{\hat{p}(s-) \sum_{k=1}^{n}\sum_{l=1}^{m_k} dM_{kl}(s)}
{\mathcal{Y}(s,\Lambda_0^\circ)},
\label{eq:mart}
\end{eqnarray}
where
 $$
\pi_r(s,\bgamma,\Lambda_0) = n^{-1} \int_s^{\tau}
 \frac{\sum_{i=1}^{n} \sum_{j=1}^{m_i}
 Q_{ijr}(\bgamma,\Lambda_0,t)
 d\tilde{N}_{ij}(t)}{\hat{p}(t)}.
 $$
Therefore, $n^{1/2}[{\bf
U}(\bgamma^\circ,\hat{\Lambda}_0(\cdot,\bgamma^\circ))-{\bf
U}(\bgamma^\circ,\Lambda_0^\circ(\cdot,\bgamma^\circ))]$
 is asymptotically mean zero multivariate normal with covariance
 matrix that can be consistently estimated by
 \begin{eqnarray*}
 G_{rl}(\hat{\bgamma}) = n^{-1} \int_0^{\tau} \pi_r(s,\hat{\bgamma},\hat{\Lambda}_0)
  \pi_l(s,\hat{\bgamma},\hat{\Lambda}_0) \{\hat{p}(s-)\}^2
 \frac{\sum_{i=1}^{n} \sum_{j=1}^{m_i} dN_{ij}(s)}
 {\{\mathcal{Y}(s,\hat{\Lambda}_0)\}^{2}
 }
 \end{eqnarray*}
 for $r,l=1,\ldots,p+1$.

\bsh
\underline{Step III}
\esh

We now examine the sum of
$\bU(\bgamma^\circ, \Lambda_0^\circ)$ and
$\bU(\bgamma^\circ,\hat{\Lambda}_0(\cdot,\bgamma^\circ))-\bU(\bgamma^\circ,
 \Lambda_0^\circ)$. From (\ref{eq:mart}), we have
 \begin{eqnarray*}
U_r(\bgamma^\circ,\hat{\Lambda}_0(\cdot,\bgamma^\circ))-U_r(\bgamma^\circ,\Lambda_0^\circ)
\approx n^{-1}\int_0^{\tau} \alpha_r(s) \sum_{k=1}^n
\sum_{l=1}^{m_k} dM_{kl}(s) =\frac{1}{n}\sum_{k=1}^n \mu_{kr},
 \end{eqnarray*}
where $\alpha_r(s)$ is the limiting value of
$\pi_r(s,\bgamma^\circ,\Lambda_0^\circ) \hat{p}(s-)
/\mathcal{Y}(s,{\Lambda}_0^\circ)$
and $\mu_{kr}$ is defined as
\begin{eqnarray*}
\mu_{kr}=\int_0^{\tau} \alpha_r(s) \sum_{l=1}^{m_k} dM_{kl}(s).
\end{eqnarray*}
Arguments in Yang and Prentice (1999, Appendix A) can
be used to show that $\hat{p}(s-)$ has a limit. Also, clearly
$\E[\mu_{kr}]=0$.

We thus have
\begin{eqnarray*}
U_r(\bgamma^\circ,
\Lambda_0^\circ)+[U_r(\bgamma^\circ,\hat{\Lambda}_0(\cdot,\bgamma^\circ))-U_r(\bgamma^\circ,
 \Lambda_0^\circ)]
 \approx
\frac{1}{n} \sum_{i=1}^n (\xi_{ir} + \mu_{ir}),
\end{eqnarray*}
which is a mean of $n$ iid random variables. Hence
$n^{1/2}\{U_r(\bgamma^\circ,
\Lambda_0^\circ)+[U_r(\bgamma^\circ,\hat{\Lambda}_0(\cdot,\bgamma^\circ))-U_r(\bgamma^\circ,
 \Lambda_0^\circ)]\}$
is asymptotically normally distributed. The covariance matrix may
be estimated by
$\hat{\bV}(\hat{\bgamma})+\hat{\bG}(\hat{\bgamma})+\hat{\bC}(\hat{\bgamma})$,
where
\begin{eqnarray*}
\hat{C}_{rl}(\hat{\bgamma})=
\frac{1}{n}
\sum_{i=1}^{n} (\xi_{ir}^{*}
\mu_{il}^{*} + \xi_{il}^{*} \mu_{ir}^{*}), \quad
r,l=1,\ldots,p+1,
\end{eqnarray*}
with
\begin{eqnarray*}
\mu_{ir}^{*}=\int_0^{\tau}
\frac{\pi_r(s,\hat{\bgamma},\hat{\Lambda}_0) \hat{p}(s-)}
{\mathcal{Y}(s,\hat{\Lambda}_0)} \sum_{j=1}^{m_i} d\hat{M}_{ij}(s)
\end{eqnarray*}
and
\begin{eqnarray*}
\hat{M}_{ij}(t)=N_{ij}(t)-\int_0^{t}\exp({\hat{\bbeta}}^T
{\bf Z}_{ij})Y_{ij}(u)
\psi_i(\hat{\bgamma},\hat{\Lambda}_0,u-)d\hat{\Lambda}_0(u).
\end{eqnarray*}

\bsh
\underline{Step IV}
\esh

First order Taylor expansion of
$\bU(\hat{\bgamma},\hat{\Lambda}_0(\cdot,\hat{\bgamma}))$ about
$\bgamma^\circ=({\bbeta^\circ}^T,\theta^\circ)^T$ gives
 \begin{eqnarray*}
\bU(\hat{\bgamma},\hat{\Lambda}_0(\cdot,\hat{\bgamma}))
=\bU(\bgamma^\circ,\hat{\Lambda}_0(\cdot,\bgamma^\circ))+ \bD(\bgamma^\circ)
(\hat{\bgamma}-\bgamma^\circ)^T + o_p(1),
 \end{eqnarray*}
where
 \begin{eqnarray*}
D_{ls}(\bgamma)=\partial
U_l(\bgamma,\hat{\Lambda}_0(\cdot,\bgamma)) / \partial \gamma_s
 \end{eqnarray*}
for $l,s=1,\ldots,p+1$, with $\gamma_{p+1}=\theta$.

For $l,s=1,\ldots,p$ we have
 \begin{eqnarray}\label{dls}
\hspace*{-2em} D_{ls}(\bgamma)&=&-n^{-1}\sum_{i=1}^{n}
 \left\{
 \frac{\phi_{2i}(\bgamma,\hat{\Lambda}_0,\tau)}{\phi_{1i}(\bgamma,\hat{\Lambda}_0,\tau)}
 \sum_{j=1}^{m_i}Z_{ijl} \frac{\partial \hat{H}_{ij}(T_{ij})}{\partial
 \beta_s} \right. \nonumber\\
& & - \left. \left[
 \frac{\phi_{3i}(\bgamma,\hat{\Lambda}_0,\tau)}{\phi_{1i}(\bgamma,\hat{\Lambda}_0,\tau)}
 -\frac{\phi^2_{2i}(\bgamma,\hat{\Lambda}_0,\tau)}{\phi^2_{1i}(\bgamma,\hat{\Lambda}_0,\tau)}
 \right]
 \sum_{j=1}^{m_i} \hat{H}_{ij}(T_{ij}) Z_{ijl}\frac{\partial \hat{H}_{i.}(\tau)}{\partial
 \beta_s}
 \right\},
 \end{eqnarray}
 \begin{eqnarray*}
\frac{\partial \hat{H}_{ij}(\tau_k)}{\partial \beta_s} =
\frac{\partial \hat{\Lambda}_0(T_{ij} \wedge \tau_k)}{\partial
\beta_s} \exp(\bbeta^T {\bf Z}_{ij}) +\hat{\Lambda}_0(T_{ij}
\wedge \tau_k)\exp(\bbeta^T {\bf Z}_{ij}) Z_{ijs}
 \end{eqnarray*}
and
\begin{eqnarray*}
\frac{\partial \Delta \hat{\Lambda}_0(\tau_k)}{\partial \beta_s}
&=& -d_k\left\{\sum_{i=1}^{n}
 \frac{\phi_{2i}(\bgamma,\hat{\Lambda}_0,\tau_{k-1})}
 {\phi_{1i}(\bgamma,\hat{\Lambda}_0,\tau_{k-1})} R_{i.}(\tau_k)
 \right\}^{-2} \nonumber \\
&& \sum_{i=1}^n \left[ \left\{
 \frac{\phi_{2i}^2(\bgamma,\hat{\Lambda}_0,\tau_{k-1})}
 {\phi_{1i}^2(\bgamma,\hat{\Lambda}_0,\tau_{k-1})}
 -
 \frac{\phi_{3i}(\bgamma,\hat{\Lambda}_0,\tau_{k-1})}
 {\phi_{1i}(\bgamma,\hat{\Lambda}_0,\tau_{k-1})}
 \right\} \frac{\partial \hat{H}_{i.}(\tau_{k-1})}{\partial \beta_s} R_{i.}(\tau_k) \right. \nonumber \\
&& + \left.
 \frac{\phi_{2i}(\bgamma,\hat{\Lambda}_0,\tau_{k-1})}
 {\phi_{1i}(\bgamma,\hat{\Lambda}_0,\tau_{k-1})}
 \sum_{j=1}^{m_i} R_{ij}(\tau_{k}) Z_{ijs}
 \right].
\end{eqnarray*}
For $l=1,\ldots,p$ we have
\begin{eqnarray}\label{dlp}
 D_{l(p+1)}(\bgamma) &=&
 -n^{-1}\sum_{i=1}^{n}\left\{
 \frac{\phi_{2i}(\bgamma,\hat{\Lambda}_0,\tau)}{\phi_{1i}(\bgamma,\hat{\Lambda}_0,\tau)}
 \sum_{j=1}^{m_i} Z_{ijl} \frac{\partial \hat{H}_{ij}(T_{ij})}{\partial \theta} \right. \nonumber\\
 &&  \left.
 + \left[ \frac{\phi_{2i}^{(\theta)}(\bgamma,\hat{\Lambda}_0,\tau)}{\phi_{1i}(\bgamma,\hat{\Lambda}_0,\tau)}
 -\frac{\phi_{2i}(\bgamma,\hat{\Lambda}_0,\tau)\phi_{1i}^{(\theta)}(\bgamma,\hat{\Lambda}_0,\tau)}
  {\phi_{1i}^2(\bgamma,\hat{\Lambda}_0,\tau)} \right. \right. \nonumber\\
 && \left. \left.
  +\left\{\frac{\phi_{2i}^2(\bgamma,\hat{\Lambda}_0,\tau)}{\phi_{1i}^2(\bgamma,\hat{\Lambda}_0,\tau)}
  - \frac{\phi_{3i}(\bgamma,\hat{\Lambda}_0,\tau)}{\phi_{1i}(\bgamma,\hat{\Lambda}_0,\tau)} \right\}
  \frac{\partial \hat{H}_{i.}(\tau)}{\partial \theta} \right]
  \sum_{j=1}^{m_i}\hat{H}_{ij}(T_{ij})Z_{ijl}
 \right\}
\end{eqnarray}
and
 \begin{eqnarray}\label{dpl}
 D_{(p+1)l}(\bgamma) &=& n^{-1}\sum_{i=1}^{n}
 \left\{
 \frac{
 \phi_{1i}^{(\theta)}(\bgamma,\hat{\Lambda}_0,\tau) \phi_{2i}(\bgamma,\hat{\Lambda}_0,\tau)}
 {\phi_{1i}^2(\bgamma,\hat{\Lambda}_0,\tau)}
-
 \frac{
 \phi_{2i}^{(\theta)}(\bgamma,\hat{\Lambda}_0,\tau)}
 {\phi_{1i}(\bgamma,\hat{\Lambda}_0,\tau)}
 \right\} \frac{\partial \hat{H}_{i.}(\tau)}{\partial \beta_l}.
 \end{eqnarray}
Finally,
 \begin{eqnarray}\label{dpp}
 D_{(p+1)(p+1)}(\bgamma)&=&n^{-1}\sum_{i=1}^{n}
 \left\{
 \frac{
 \phi_{1i}^{(\theta,\theta)}(\bgamma,\hat{\Lambda}_0,\tau)
 }{\phi_{1i}(\bgamma,\hat{\Lambda}_0,\tau)}
  -
\left[ \frac{
 \phi_{1i}^{(\theta)}(\bgamma,\hat{\Lambda}_0,\tau)
 }{\phi_{1i}(\bgamma,\hat{\Lambda}_0)}\right]^2 \right. \nonumber \\
&& \left. +\left[ \frac{
\phi_{1i}^{(\theta)}(\bgamma,\hat{\Lambda}_0,\tau)
\phi_{2i}(\bgamma,\hat{\Lambda}_0,\tau)}{
\phi_{1i}^{2}(\bgamma,\hat{\Lambda}_0,\tau)}
 -\frac{\phi_{2i}^{(\theta)}(\bgamma,\hat{\Lambda}_0,\tau)}{\phi_{1i}(\bgamma,\hat{\Lambda}_0,\tau)} \right]
 \frac{\partial \hat{H}_{i.}(\tau)}{\partial \theta}
 \right\}
 \end{eqnarray}
where
 $$
\phi_{1i}^{(\theta,\theta)}(\bgamma,\hat{\Lambda}_0,\tau) =\int
w^{N_{i.}(\tau)}\exp\{-w \hat{H}_{i.}(\tau)\} \frac{d^2 f(w)}
 {d \theta^2} dw,
 $$
 \begin{eqnarray*}
\frac{\partial \hat{H}_{ij}(\tau_k)}{\partial \theta} =
\frac{\partial \hat{\Lambda}_0(T_{ij} \wedge \tau_{k})}{\partial
\theta} \exp(\bbeta^T {\bf Z}_{ij}),
 \end{eqnarray*}
and
\begin{eqnarray*}
\frac{\partial \Delta \hat{\Lambda}_0(\tau_k)}{\partial \theta}
&=& -d_k\left\{\sum_{i=1}^{n}
 \frac{\phi_{2i}(\bgamma,\hat{\Lambda}_0,\tau_{k-1})}
 {\phi_{1i}(\bgamma,\hat{\Lambda}_0,\tau_{k-1})} R_{i.}(\tau_k)
 \right\}^{-2} \nonumber \\
&& \sum_{i=1}^n
 R_{i.}(\tau_k)
 \left[
 \frac{\phi_{2i}^{(\theta)}(\bgamma,\hat{\Lambda}_0,\tau_{k-1})}
 {\phi_{1i}(\bgamma,\hat{\Lambda}_0,\tau_{k-1})}
 -
 \frac{\phi_{2i}(\bgamma,\hat{\Lambda}_0,\tau_{k-1})\phi_{1i}^{(\theta)}(\bgamma,\hat{\Lambda}_0,\tau_{k-1})}
 {\phi_{1i}^2(\bgamma,\hat{\Lambda}_0,\tau_{k-1})}
  \right.  \nonumber \\
 && \left.
 + \frac{\partial \hat{H}_{i.}(\tau_{k-1})}{\partial \theta}
 \left\{
  \frac{\phi_{2i}^2(\bgamma,\hat{\Lambda}_0,\tau_{k-1})}{\phi_{1i}^2(\bgamma,\hat{\Lambda}_0,\tau_{k-1})}
  - \frac{\phi_{3i}(\bgamma,\hat{\Lambda}_0,\tau_{k-1})}{\phi_{1i}(\bgamma,\hat{\Lambda}_0,\tau_{k-1})}
  \right\} \right].
\end{eqnarray*}

\bsh
\underline{Step V}
\esh

Combining the results above we get that
$n^{1/2}(\hat{\bgamma}-\bgamma^\circ)$ is asymptotically zero-mean
normally distributed with a covariance matrix that can be
consistently estimated by
 \begin{eqnarray*}
  \hat{\bD}^{-1}(\hat{\bgamma})\{\hat{\bV}(\hat{\bgamma})+\hat{\bG}(\hat{\bgamma})+\hat{\bC}(\hat{\bgamma})\}
  \hat{\bD}^{{-1}}(\hat{\bgamma})^T.
 \end{eqnarray*}

\subsection{Definition and behavior of $h_i(r,s)$}

The quantity $h_i(r,s)$ appearing in (\ref{eq:yapp}) is given by
\begin{eqnarray*}
h_i(r,s)&=&\frac{2R_{i.}(s)\eta_{1i}(r,s)}{\{{\mathcal
Y}(s,\Lambda_0^\circ+r\Delta)\}^3}\frac{1}{n}\sum_{l=1}^{n}R_{l.}(s)\eta_{1l}(r,s)\sum_{j=1}^{m_i}
\exp(\bbeta^T {\bf Z}_{lj})\Delta(T_{lj}\wedge s) \nonumber \\
 && - \frac{R_{i.}(s)\eta_{2i}(r,s)}{\{{\mathcal
Y}(s,\Lambda_0^\circ+r\Delta)\}^2}\sum_{j=1}^{m_i} \exp(\bbeta^T {\bf
Z}_{ij})\Delta(T_{ij}\wedge s)
\end{eqnarray*}
where $\Delta(T_{ij}\wedge s)=\hat{\Lambda}_0(T_{ij}\wedge
s)-\Lambda_0^o(T_{ij}\wedge s)$ and
 $$
\eta_{2i}(r,s)=2\left\{\frac{\phi_{2i}(\bgamma^\circ,{\Lambda}_0^\circ+r\Delta,s)}
{\phi_{1i}(\bgamma^\circ,{\Lambda}_0^\circ+r\Delta,s)}\right\}^3
+\frac{\phi_{4i}(\bgamma^\circ,{\Lambda}_0^\circ+r\Delta,s)}
{\phi_{1i}(\bgamma^\circ,{\Lambda}_0^\circ+r\Delta,s)}
-3\frac{\phi_{2i}(\bgamma^\circ,{\Lambda}_0^\circ+r\Delta,s)\phi_{3i}(\bgamma^\circ,{\Lambda}_0^\circ+r\Delta,s)}
{\left\{\phi_{1i}(\bgamma^\circ,{\Lambda}_0^\circ+r\Delta,s)\right\}^2}.
 $$
For all $i=1,\ldots, n$ and $s \in [0,\tau]$,
we have $0 \leq R_{i.}(s) \leq m \nu$,
where $\nu$ is as in (\ref{eb}).
Moreover, for $k=1,\ldots,4$, we have
 $$
 \E[W_i^{r_{min}+(k-1)}\exp\{-W_i m e^{\sbeta^T Z} \Lambda_0^\circ(\tau)\}]
 \leq \phi_{ki}(\bgamma^\circ,{\Lambda}_0^\circ,s) \leq
 \E[W_i^{r_{max}+(k-1)}]
  $$
where $r_{\max}=\arg \max_{1\leq r \leq m} \E (W_i^r)$,
$r_{\min}=\arg \min_{1\leq r \leq m} \E (W_i^r)$. Hence,
$\eta_{1i}$ and $\eta_{2i}$ are bounded. In addition, the the
proof of Lemma 2 show that ${\mathcal Y}(s,\Lambda^\circ
+r\Delta)$ is uniformly bounded away from zero for $n$
sufficiently large. Finally, in the consistency proof we obtained
$\| \Delta \| = o(1)$. Therefore $h_{i}(r,s)$ is $o(1)$ uniformly
in $r$ and $s$.

\section{References}
\begin{description}
\item
{\sc Aalen, O. O.} (1976). Nonparametric inference in connection
with multiple decrement models. {\em Scand. J. Statist.} {\bf 3},
15-27.
\item
{\sc Abramowitz, M. and Stegun, I. A. (Eds.)} (1972). {\em
Handbook of Mathematical Functions with Formulas, Graphs, and
Mathematical Tables, 9th printing} New York: Dover.
\item
{\sc Andersen, P. K., Borgan, O, Gill, R. D. and Keiding, N.}
(1993). Statistical models based on counting processes. Berlin:
Springer-Verlag.
\item
{\sc Andersen, P. K. and Gill, R. D.} (1982). Cox's regression
model for counting processes: A large sample study. {\em  Ann.
Statist.} {\bf 10}, 1100-1120.
\item
{\sc Andersen, P. K., Klein, J. P., Knudsen, K. M. and Palacios,
R. T.} (1997). Estimation of variance in Cox's regression model
with shared gamma frailty. {\em Biometrics} {\bf 53}, 1475-1484.
\item
{\sc Breslow, N.} (1974). Covariance analysis of censored survival
data. {\em Biometrics,} {\bf 30}, 89-99.
\item
{\sc Fine, J. P., Glidden D. V. and Lee, K.} (2003). A simple
estimator for a shared frailty regression model. {\em J. R.
Statist. Soc.} {\bf B 65}, 317-329.
\item
{\sc Cox, D. R.} (1972). Regression models and life tables (with
discussion). {\em J. R. Statist. Soc.} {\bf B 34}, 187-220.
\item
{\sc Foutz, R. V.} (1977). On the unique consistent solution to
the likelihood equation. {\em J. Amer. Statist. Ass.} {\bf 72},
147-148.
\item
{\sc Gill, R. D.} (1985). Discussion of the paper by D. Clayton
and J. Cuzick. {\em J. R. Statist. Soc.} {\bf A 148}, 108-109.
\item
{\sc Gill, R. D.} (1989). Non- and semi-parametric maximum
likelihood estimators and the Von Mises method (Part 1). {\em
Scand. J. Statist.} {\bf 16}, 97-128.
\item
{\sc Gill, R. D.} (1992). Marginal partial likelihood. {\em Scand.
J. Statist.} {\bf 79}, 133-137.
\item
{\sc Hartman, P.} (1973). {\em Ordinary Differential Equations,}
2nd ed. (reprinted, 1982), Boston: Birkhauser.
\item
{\sc Henderson, R. and Oman, P.} (1999). Effect of frailty on
marginal regression estimates in survival analysis. {\em J. R.
Statist. Soc.} {\bf B 61}, 367-379.
\item
{\sc Hougaard, P.} (1986). Survival models for heterogeneous
populations derived from stable distributions. {\em Biometrika}
{\bf 73}, 387-396.
\item
{\sc Hougaard, P.} (2000). {\em Analysis of Multivariate Survival
data}. New York: Springer.
\item
{\sc Hsu, L., Chen, L., Gorfine, M. and Malone, K.} (2004).
Semiparametric estimation of marginal hazard function from
case-control family studies. {\em Biometrics} {\bf 60}, 936-944.
\item
{\sc Klein, J. P.} (1992). Semiparametric estimation of random
effects using the Cox model based on the EM Algorithm. {\em
Biometrics} {\bf 48}, 795-806.
\item
{\sc Louis, T. A.} (1982). Finding the observed information matrix
when using the EM algorithm. {\em J. R. Statis. Soc.} {\bf B 44},
226-233.
\item
{\sc McGlichrist, C. A.} (1993). REML estimation for survival
models with frailty. {\em Biometrics} {\bf 49}, 221-225.
\item
{\sc Murphy, S. A.} (1994). Consistency in a proportional hazards
model incorporating a random effect. {\em Ann. Statist.} {\bf 22},
712-731.
\item
{\sc Murphy, S. A.} (1995). Asymptotic theory for the frailty
model. {\em Ann. Statist.} {\bf 23}, 182-198.
\item
{\sc Nielsen, G. G., Gill, R. D., Andersen, P. K. and Sorensen, T.
I.} (1992). A counting process approach to maximum likelihood
estimation of frailty models. {\em Scand. J. Statist.} {\bf 19},
25-43.
\item
{\sc Parner, E.} (1998). Asymptotic theory for the correlated
gamma-frailty model. {\em Ann. Statist.} {\bf 26}, 183-214.
\item
{\sc Ripatti, S. and Palmgren J.} (2000). Estimation of
multivariate frailty models using penalized partial likelihood.
{\em Biometrics} {\bf 56}, 1016-1022.
\item
{\sc Shih, J. H. and Chatterjee, N.} (2002). Analysis of survival
data from case-control family studies. {\em Biometrics} {\bf 58},
502-509.
\item
{\sc Vaida, F. and Xu, R. H.} (2000). Proportional hazards model
with random effects. {\em Stat. in Med.} {\bf 19}, 3309-3324.
\item
{\sc Yang, S. and Prentice, R. L.} (1999). Semiparametric
inference in the proportional odds regression model. {\em J. Amer.
Statist. Ass.} {\bf 94}, 125-136.
\item
{\sc Zucker, D. M.} (2005). A pseudo partial likelihood method for
semi-parametric survival regression with covariate errors. {\em to
appear in J. Amer. Statist. Ass.}.
\end{description}

\newpage
\begin{center}
\renewcommand{\baselinestretch}{1.2}
\begin{table}
\caption{ {\em Simulation results for the gamma frailty model with
single normal covariate; $n=300$ and family size equals 2. $Z \sim
N(0,1)$; $\beta^\circ=\log(2)$ or $\log(3)$; $\theta^\circ=2$ }}
\begin{center}
\begin{tabular}{lcccc}
\hline \hline
  & \multicolumn{2}{c}{$\hat{\beta}$} & \multicolumn{2}{c}{$\hat{\theta}$}\\
  & \multicolumn{2}{c}{\rule{5.0cm}{0.05mm}} &
\multicolumn{2}{c}{\rule{5.0cm}{0.05mm}} \\
 & Our approach & EM algorithm &Our approach & EM algorithm \\
 \hline
  &  \multicolumn{4}{c}{$\beta=ln(2)$; 35\% censoring} \\
 Empirical mean & 0.692 & 0.689 & 1.978 & 1.969 \\
 Empirical SD & 0.248 & 0.253 & 0.268 & 0.308\\
 Mean estimated SD & 0.242 & - & 0.242 & -\\
 $95\%$ Wald-type CI & 95.6 & - & 96.3 & - \\
 Correlation & \multicolumn{2}{c}{0.952} & \multicolumn{2}{c}{0.989}\\
  &  \multicolumn{4}{c}{$\beta=ln(2)$; 85\% censoring} \\
 Empirical mean & 0.699 & 0.693 & 1.942 & 1.942 \\
 Empirical SD & 0.479 & 0.481 & 0.897 & 0.936 \\
 Mean estimated SD & 0.442 & - & 0.919 & - \\
 $95\%$ Wald-type CI & 96.6 & - & 95.0 & - \\
 Correlation & \multicolumn{2}{c}{0.952} & \multicolumn{2}{c}{0.989}\\
 &  \multicolumn{4}{c}{$\beta=ln(3)$; 30\% censoring} \\
Empirical mean & 1.102 & 1.078 & 1.985 & 1.961 \\
 Empirical SD & 0.255 & 0.266 & 0.265 & 0.259\\
 Mean estimated SD & 0.231 & - & 0.279 & -\\
 $95\%$ Wald-type CI & 96.9 & - & 96.1 & - \\
 Correlation & \multicolumn{2}{c}{0.951} & \multicolumn{2}{c}{0.982}\\
  &  \multicolumn{4}{c}{$\beta=ln(3)$; 80\% censoring} \\
 Empirical mean & 1.099 & 1.088 & 1.921 & 1.870 \\
 Empirical SD & 0.465 & 0.466 & 0.800 & 0.810 \\
 Mean estimated SD & 0.443 & - & 0.797 & - \\
 $95\%$ Wald-type CI & 94.2 & - & 96.3 & - \\
 Correlation & \multicolumn{2}{c}{0.957} & \multicolumn{2}{c}{0.993}\\
\hline
\end{tabular}
\end{center}
\end{table}
\end{center}

\end{document}